\documentclass[final,1p,times,authoryear]{elsarticle}

\usepackage{amsmath}
\usepackage{amssymb}

\newcommand{\refeqn}[1]{{\rm (\ref{#1})}}

\journal{Journal of Computational Physics}

\begin{document}

\begin{frontmatter}



\title{A primal-dual mimetic finite element scheme for the rotating
shallow water equations on polygonal spherical meshes}


\author[au_jt]{John Thuburn}
\address[au_jt]{College of Engineering, Mathematics and Physical Sciences, University of Exeter,
Exeter, UK}
\ead{j.thuburn@exeter.ac.uk}

\author[au_cjc]{Colin J Cotter}
\address[au_cjc]{Department of Mathematics, Imperial College, London, UK}

\begin{abstract}
A new numerical method is presented for solving the rotating shallow water
equations on a rotating sphere using quasi-uniform polygonal meshes.
The method uses special families of finite element function spaces to
mimic key mathematical properties of the continuous equations and thereby
capture several desirable physical properties related to balance and
conservation. The method relies on two novel features. The first is the use
of {\it compound finite elements} to provide suitable finite element spaces
on general polygonal meshes. The second is the use of
{\it dual finite element spaces} on the dual of the original mesh, along with
suitably defined discrete Hodge star operators to map between the primal
and dual meshes, enabling the use of a finite volume scheme on the dual
mesh to compute potential vorticity fluxes.
The resulting method has the same mimetic properties as a finite volume
method presented previously, but is more accurate on a number of standard
test cases.

\end{abstract}


\end{frontmatter}


\section{Introduction}
\label{sec_introduction}

In order to exploit the new generation of massively parallel supercomputers
that are becoming available, weather and climate
models will require good parallel scalability. This requirement has driven the
development of numerical methods that do not depend on the orthogonal
coordinate system and quadrilateral structure of the longitude-latitude grid,
whose polar resolution clustering is predicted to lead to a scalability
bottleneck.
A significant challenge is to obtain good scalability without
sacrificing accuracy; in particular conservation, balance, and wave propagation
are important for accurate modelling of the atmosphere \citep{staniforth2012}.

Building on earlier work \citep{ringler2010,thuburn2012}, \citet{thuburn2014}
presented a finite volume scheme for the shallow water equations on polygonal
meshes. They start from the continuous shallow water equations in the
so-called vector invariant form:
\begin{eqnarray}
\phi_t + \nabla \cdot \mathbf{f}  & = & 0, \label{ctsmass} \\
\mathbf{u}_t + \mathbf{q}^{\perp} + \nabla (\phi_{\mathrm{T}} + k) & = & 0,
\label{ctsvel}
\end{eqnarray}
where $\phi$, the geopotential, is equal to the fluid depth times the
gravitational acceleration, $\phi_{\mathrm{T}} = \phi +
\phi_{\mathrm{orog}}$ is the total geopotential at the fluid's upper
surface including the contribution from orography, $\mathbf{u}$ is the
velocity, $\mathbf{f} = \mathbf{u} \phi$ is the mass flux,
and $k = |\mathbf{u}|^2/2$.  The $\perp$ symbol is defined by
$\mathbf{u}^{\perp} = \mathbf{k} \times \mathbf{u}$ where $\mathbf{k}$ is the unit
vertical vector.  Finally,  $\pi = \zeta/\phi$ is the potential vorticity
(PV), where $\zeta = f + \xi$ is the absolute
vorticity, with $f$ the Coriolis parameter and $\xi = \mathbf{k} \cdot \nabla \times \mathbf{u}$
the relative vorticity, and $\mathbf{q} = \mathbf{f} \pi$ is the PV flux.
By the use of a C-grid placement of prognostic variables, and by ensuring
that the numerical method mimics key mathematical properties of the continuous
governing equations (hence the term `mimetic'), the scheme was designed to have
good conservation and balance properties. These good properties were verified
in numerical tests on hexagonal and cubed sphere spherical meshes.
However, their scheme has a number of drawbacks. Most seriously, the Coriolis
operator, whose discrete form is essential to obtaining good geostrophic balance,
is numerically inconsistent and fails to converge in the $L_\infty$ norm
\citep{weller2014,thuburn2014}.
Also, although the gradient and divergence operators are consistent,
their combination to form the discrete Laplacian operator also fails to converge
in the $L_\infty$ norm in some cases. These inaccuracies are clearly visible in
idealized convergence tests, and give rise to marked `grid imprinting'
for initially symmetrical flows.
Although they are less conspicuous in more complex flows, they are clearly
undesirable.

\citet{cotter2012} \citep[see also][]{mcrae2014,cotter2014} showed that the same
mimetic properties can be obtained using a certain class of mixed finite
element method. The mimetic properties follow from the choice of an appropriate
hierarchy of function spaces for the prognostic and diagnostic variables
(e.g.~section~\ref{sec_compound} below), which also provides a finite element analogue
of the C-grid placement of variables, or a higher-order generalization.
(The use of such a hierarchy goes by various names in the literature, including
`mimetic finite elements', `compatible finite elements', and `finite element
exterior calculus'; see \citet{cotter2014} for a discussion of the shallow water
equation case in the language of exterior calculus.)
Importantly, the resulting schemes are numerically consistent.

While the mimetic finite element approach appears very attractive, it is not
yet clear which particular choice of mesh and function spaces is most suitable.
Standard finite element methods use
triangular or quadrilateral elements. For the lowest-order mimetic finite element
scheme on triangles, the dispersion relation for the linearized shallow water
equations suffers from extra branches of inertio-gravity waves, which are
badly behaved numerical artefacts \citep{leroux2007}, analogous to the problem
that occurs on the triangular C-grid \citep{danilov2010}. Higher-order finite
element methods also typically exhibit anomalous features in their wave dispersion
relations, such as extra branches, frequency gaps, or zero group velocity modes.
Some progress has been made in reducing these
problems, at least on quadrilateral meshes, through the inclusion of dissipation
or modification of the mass matrix \citep[e.g.][]{melvin2013,ullrich2013},
though the remedies are somewhat heuristic except in the most idealized cases.
Finally, coupling to subgrid models of physical processes such as
cumulus convection or cloud microphysics may be less straightforward with
higher-order elements (P.\ Lauritzen, pers.\ comm.).
These factors suggest that it may still be worthwhile investigating
lowest-order schemes on quadrilateral and hexagonal meshes.

The above arguments raise two related questions.
Can the mimetic finite element method inspire a development to
fix the inconsistency of the mimetic finite volume method? 
Alternatively, can the mimetic finite element method at lowest order
be adapted to work on polygonal meshes such as hexagons?
Below we answer the second question by showing that the mimetic finite
element method can indeed be adapted. 
In fact, from a certain viewpoint the mimetic finite volume and mimetic
finite element schemes have very similar mathematical structure. 
The notation below is chosen to emphasize this
similarity\footnote{Readers wishing to compare the two formulations
should note that a different sign convention is used for the expansion
coefficients of $\mathbf{k}\times$ any vector, such as $U^\perp$
in \refeqn{uperp}.}.
Moreover,
the similarity is sufficiently strong that much of the code of the mimetic
finite volume model of \citet{thuburn2014} could be re-used in the model
presented below. This, in turn, facilitates the cleanest possible comparison
of the two approaches.

The adaptation of the mimetic finite element method employs two novel
features. The first is the definition of a suitable hierarchy of
finite element
function spaces on polygonal meshes. This is achieved by defining compound
elements built out of triangular subelements, and is described in
section~\ref{sec_compound}. The second ingredient is the introduction of
a dual family of function spaces that are defined on the dual of the
original mesh. This permits the definition of a spatially averaged mass field
that lives in the same function space as the vorticity and potential vorticity
fields; this, in turn, enables the use of an accurate finite volume scheme on the dual
mesh for advection of potential vorticity, and keeps the formulation of the
finite element model as close as possible to that of the finite volume model.

\section{Meshes and dual meshes}
\label{sec_meshes}

\begin{figure}
\includegraphics[width=60mm]{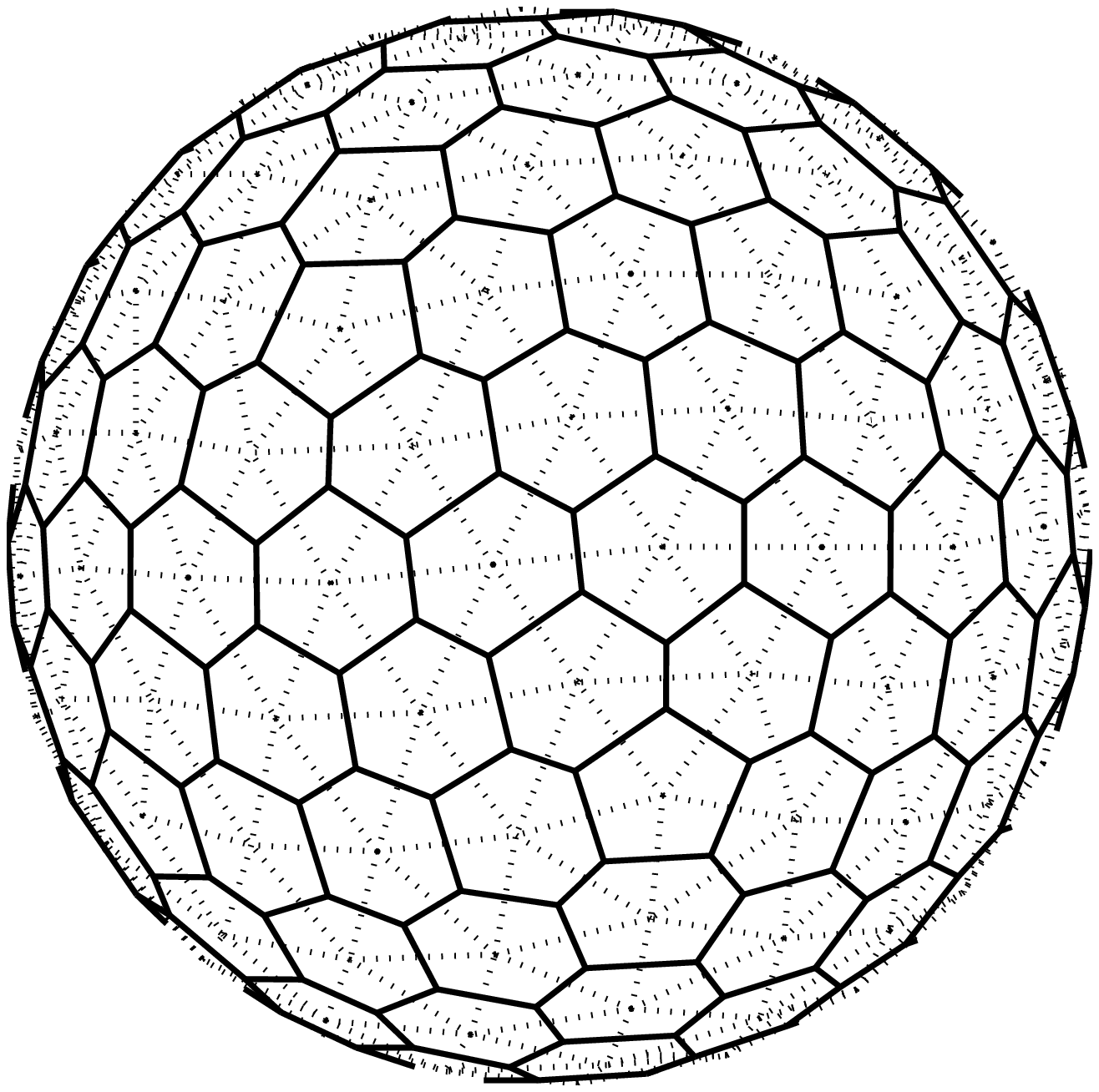}
\includegraphics[width=60mm]{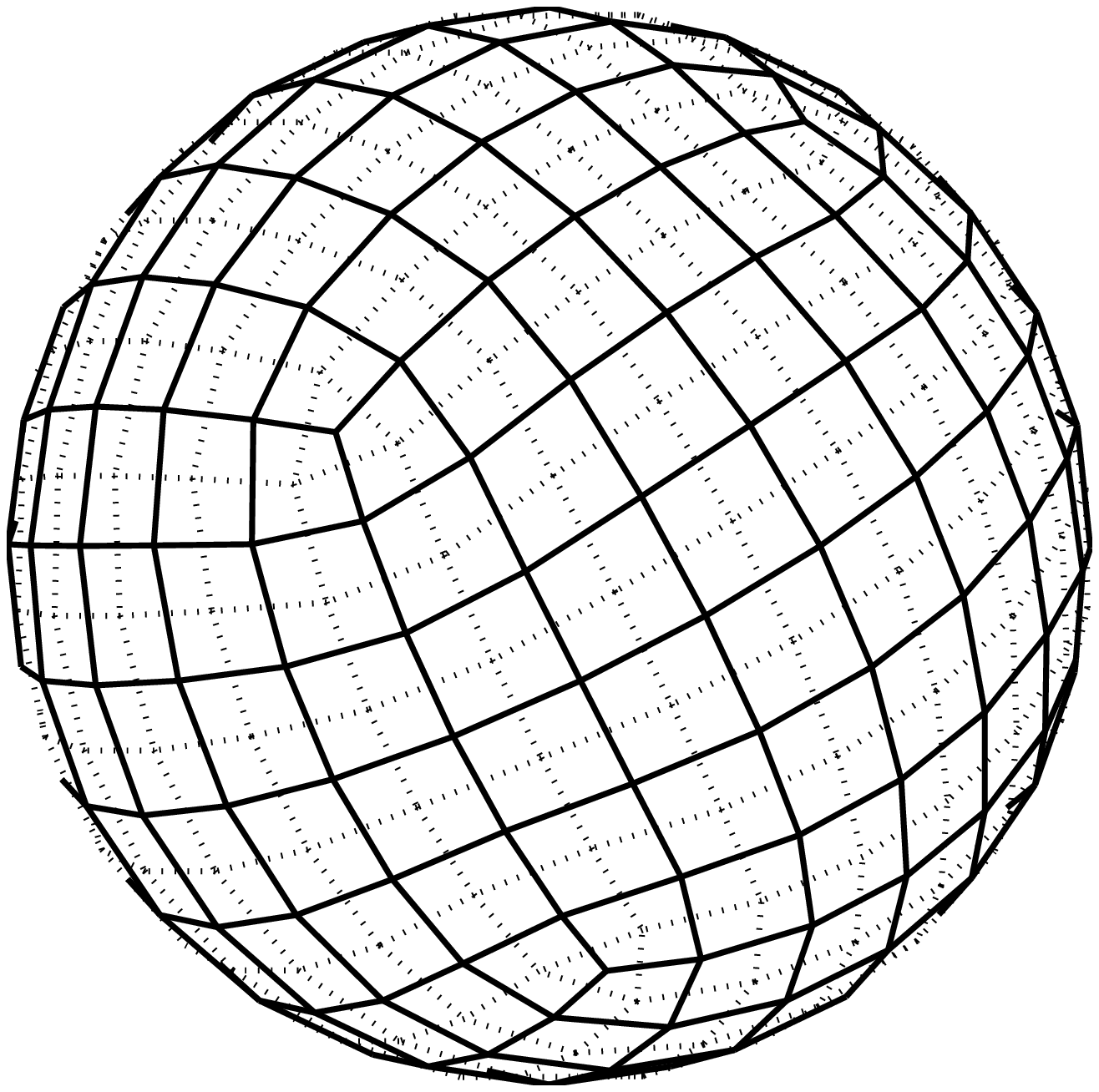}
\caption{Left: a~hexagonal--icosahedral mesh with 162 cells and
  642 degrees of freedom. Right: a~cubed-sphere mesh with 216 cells
  and 648 degrees of freedom.  Continuous lines are primal mesh edges,
  dotted lines are dual mesh edges.}
\label{fig_grids}
\end{figure}

The scheme described here is suitable for arbitrary two-dimensional polygonal meshes
on flat domains or, as used here, curved surfaces approximated by planar facets.
Two particular meshes are used to obtain the results in section~\ref{sec_results},
namely the same variants of the hexagonal-icosahedral mesh and the cubed sphere mesh
used by \citet{thuburn2014}, in order to facilitate comparison with their results.
Coarse-resolutions versions are shown in Fig.~\ref{fig_grids}.

Any polygonal mesh has a corresponding dual mesh. (We will refer to the original mesh
as the `primal' mesh where necessary to distinguish it from the dual.) Each primal
cell contains one dual vertex; each dual cell contains one primal vertex; each
primal edge corresponds to one dual edge and these usually cross each other.
Figure~\ref{fig_grids} shows both primal and dual edges for the two meshes.

\section{Function spaces and compound finite elements}
\label{sec_compound}

\begin{figure}
\includegraphics[width=150mm]{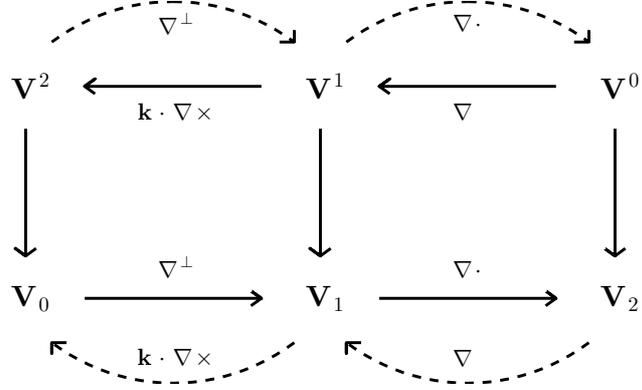}
\caption{
Schematic showing the function spaces used in the scheme
and the relationships between them. Primal function spaces are on the
bottom row and dual function spaces are on the top row.}
\label{fig_spaces}
\end{figure}

The mimetic properties of the scheme arise from the relationships between the
finite element function spaces. Three function spaces are used on the primal mesh
($\mathbf{V}_0$, $\mathbf{V}_1$, and $\mathbf{V}_2$), and three on the dual mesh
($\mathbf{V}^2$, $\mathbf{V}^1$, and $\mathbf{V}^0$).
Figure~\ref{fig_spaces} indicates that $\nabla^\perp$ (i.e.\ $\mathbf{k} \times \nabla$)
maps from
$\mathbf{V}_0$ to $\mathbf{V}_1$ and $\nabla \cdot$ maps
from $\mathbf{V}_1$ to $\mathbf{V}_2$\footnote{Note
$\nabla^\perp$ and $\mathbf{k} \cdot \nabla \times$
(like $\nabla$ and $\nabla \cdot$) can both be
defined as intrinsic operations on a curved surface, without reference to
$\mathbf{k}$ or a third dimension.}.
More precisely,
the primal function spaces satisfy the following properties.

\vspace{3mm}
\noindent
\textbf{Property List 1}
\begin{itemize}
\item
$\mathbf{u} \in \mathbf{V}_1 \ \ \Rightarrow \ \ 
\nabla \cdot \mathbf{u} \in \mathbf{V}_2$.

\item
$\phi \in \mathbf{V}_2$ with 
$\int \phi \, dA = 0 \ \ \Rightarrow \ \ \exists \,
\mathbf{u} \in \mathbf{V}_1 \mathrm{\ st\ } \nabla \cdot \mathbf{u} = \phi$.

\item
$\psi \in \mathbf{V}_0 \ \ \Rightarrow \ \ \nabla^\perp \psi \in \mathbf{V}_1$.

\item
$\forall \psi \in \mathbf{V}_0 , \ \nabla \cdot \nabla^\perp \psi = 0$, and
$\forall \mathbf{u} \in \mathbf{V}_1 \ \mathrm{st} \ \nabla \cdot \mathbf{u} = 0, \ 
\exists \, \psi \in \mathbf{V}_0 \ \mathrm{st} \ \mathbf{u} = \nabla^\perp \psi$.
That is, $\nabla^\perp$ maps onto the kernel of $\nabla \cdot$\ .

\end{itemize}
The second condition assumes spherical geometry so that there are no lateral
boundaries. The same assumption will be made throughout this paper;
in particular, no boundary terms will arise when integrating by
parts\footnote{The most general form of the Helmholtz-Hodge decomposition of a vector
field $\mathbf{u}$ in 2D is
\begin{displaymath}
\mathbf{u} = \nabla \phi + \nabla^\perp \psi + \mathbf{h}, 
\end{displaymath}
where $\phi$ is a potential, $\psi$ is a stream function, and
$\mathbf{h}$ is a harmonic vector field, i.e.\ one satisfying
\begin{displaymath}
\nabla \nabla \cdot \mathbf{h} +
\nabla^\perp \mathbf{k} \cdot \nabla \times \mathbf{h} =
\mathbf{0}.
\end{displaymath}
The fourth condition in Property List~1 implies that all nondivergent fields
$\mathbf{u}$ can be written as $\nabla^\perp \psi$, which rules out
the possibility of harmonic vector fields. This is appropriate for spherical
geometry, since there exist no non-zero harmonic vector fields on the sphere.
However, for a doubly period plane, for example, for which a constant
vector field is harmonic, we would have to extend the fourth condition
to allow for harmonic vector fields. This issue does not affect any
of the discussion below except for the discrete Helmholtz
decomposition (section~\ref{deriv2}), which would only need to be extended in the obvious
way to allow for harmonic vector fields.}.

In a similar way, Fig.~\ref{fig_spaces} indicates that $\nabla$ maps from
$\mathbf{V}^0$ to $\mathbf{V}^1$ and $\mathbf{k} \cdot \nabla \times$
maps from $\mathbf{V}^1$ to $\mathbf{V}^2$. More precisely, the dual function spaces
satisfy the following properties.

\vspace{3mm}
\noindent
\textbf{Property List 2}
\begin{itemize}
\item
$\hat{\mathbf{u}} \in \mathbf{V}^1 \ \ \Rightarrow \ \ 
\mathbf{k} \cdot \nabla \times \hat{\mathbf{u}} \in \mathbf{V}^2$.

\item
$\hat{\xi} \in \mathbf{V}^2$ with 
$\int \hat{\xi} \, dA = 0 \ \ \Rightarrow \ \ \exists \,
\hat{\mathbf{u}} \in \mathbf{V}^1 \mathrm{\ st\ }
\mathbf{k} \cdot \nabla \times \hat{\mathbf{u}} = \hat{\xi}$.

\item
$\hat{\chi} \in \mathbf{V}^0 \ \ \Rightarrow \ \ \nabla \hat \chi \in \mathbf{V}^1$.

\item
$\forall \hat{\chi} \in \mathbf{V}^0, \ \mathbf{k} \cdot \nabla \times \nabla \hat{\chi} = 0$,
and
$\forall \hat{\mathbf{u}} \in \mathbf{V}^1 \ \mathrm{st} \ \mathbf{k} \cdot \nabla \times
\hat{\mathbf{u}} = 0, \ \exists \, \hat{\chi} \in \mathbf{V}^0 \ \mathrm{st} \ 
\hat{\mathbf{u}} = \nabla \hat{\chi}$.
That is, $\nabla$ maps onto the kernel of $\mathbf{k} \cdot \nabla \times$\ .

\end{itemize}

As noted earlier, standard finite element schemes in two dimensions typically use
triangular or quadrilateral elements. Several families of mixed finite elements
that satisfy Property List~1 on such meshes are known.
However, in order to apply our scheme on more general
polygonal meshes we will need to define families of mixed finite elements
satisfying Property List~1 on those meshes. One way to do this is
to use \textit{compound elements}. Any polygonal element can be subdivided into a
number of triangular subelements. A basis function on the polygonal element can
then be defined as a suitable linear combination of basis functions on the
subelements. The allowed linear combinations are determined by the requirement
to satisfy Property List~1 or~2; see below.

The desire to use a dual mesh increases the need for finite element spaces
on polygons, and hence for compound elements. Only in special cases
(such as the cubed sphere, Fig.~\ref{fig_grids}) can both the primal and dual
meshes be built of triangles and quadrilaterals; other cases require higher
degree polygons for either the primal or dual mesh (or both).

For a triangular primal mesh, \citet{buffa2007} describe a scheme for the
construction of a dual hierarchy of function spaces. The dual mesh elements
are compound elements, similar, though not identical, to those used here.
However, their scheme is limited to the case of a triangular primal mesh and a
barycentric refinement for the construction of the dual.

In a complementary study, \citet{christiansen2008} describes how finite element
basis functions satifying Property List~1 may be constructed on arbitrary
polygonal elements, without the need to divide into subelements, through a process
of \textit{harmonic extension}. For example, let $\gamma_j$ be a basis function for
$\mathbf{V}_0$ associated with primal vertex~$j$. Define $\gamma_j$ to equal~$1$ at
vertex~$j$ and zero at all other vertices. Next extend $\gamma_j$ harmonically
along primal mesh edges; that is, its second derivative should vanish so that
its gradient is constant along each edge. Then extend $\gamma_j$ harmonically
into the interior of each element; that is, solve
\begin{equation}
\label{harmonicV0}
\nabla^2 \gamma_j = 0
\end{equation}
subject to the Dirichlet boundary conditions given by the known values
of $\gamma_j$ on element edges.
In a similar way, let $\mathbf{v}_e$ be a basis function for $\mathbf{V}_1$ associated
with edge~$e$. Define the normal component of $\mathbf{v}_e$ to be a
nonzero constant
along edge~$e$ (some arbitrary sign convention must be chosen to define the
positive direction) and zero at all other edges. Then extend $\mathbf{v}_e$
harmonically into the interior of each element; that is, solve
\begin{equation}
\label{harmonicV1a}
\nabla \left( \nabla \cdot \mathbf{v}_e \right) = 0
\end{equation}
and
\begin{equation}
\label{harmonicV1b}
\nabla^\perp \left( \mathbf{k} \cdot \nabla \times \mathbf{v}_e \right) = 0
\end{equation}
subject to the known values of the normal component at element edges,
for example by writing $\mathbf{v}_e = \nabla \phi + \nabla^\perp \psi$,
implying $\nabla^2 \phi = c_1$ and $\nabla^2 \psi = c_2$ for constants
$c_1$ and $c_2$. The boundary conditions determine the value of $c_1$,
but not $c_2$. 
However, condition~\refeqn{harmonicV0}
along with the fourth property in List~1 implies that we must choose
$c_2 = 0$, so that \refeqn{harmonicV1b} reduces to
\begin{equation}
\label{harmonicV1c}
\mathbf{k} \cdot \nabla \times \mathbf{v}_e = 0 .
\end{equation}
For the last function space $\mathbf{V}_2$ the basis function associated with
cell~$i$ is defined to be a nonzero constant in cell~$i$ and zero in all other cells.
It may then be verified that the properties in List~1 do indeed hold for
the spaces spanned by these basis functions.

Although the harmonic extension approach provides a general method for
constructing the lowest order mimetic finite element spaces on
polygonal meshes, its drawback is that, except for the simplest
element shapes, the basis functions cannot be found analytically. Even if they
are found numerically, the inner products required for the finite element method
cannot be computed exactly, either analytically or by numerical quadrature.

Here we take inspiration from both \citet{buffa2007} and \citet{christiansen2008}
to construct spaces of compound finite elements for arbitrary polygonal
primal and dual meshes, by a process that might be called
\textit{discrete harmonic extension}.
For the function spaces on the primal mesh, in effect, we solve a
finite element discretization of \refeqn{harmonicV0},
\refeqn{harmonicV1a}, and \refeqn{harmonicV1c} on the mesh of triangular
subelements in order to construct the compound basis elements for the original
polygonal mesh. For this discretization we use the lowest order mimetic
finite element spaces on the triangular subelements, in which
$\mathbf{V}_0$ comprises continuous piecewise linear elements,
$\mathbf{V}_1$ comprises the lowest order Raviart-Thomas elements,
and $\mathbf{V}_2$ comprises piecewise constant elements;
$\text{P1-RT0-P0}^\text{DG}$ in standard shorthand.
Although only discrete versions of \refeqn{harmonicV0},
\refeqn{harmonicV1a}, and \refeqn{harmonicV1c} are solved,
it may be verified that the properties in List~1 hold exactly.
The basis functions on the triangular subelements are known
analytically, and the compound elements are linear combinations of these;
therefore, integrals of products of basis functions, for example to compute
entries of a mass matrix, can all be computed exactly.

The resulting compound elements provide a generalization to polygonal
meshes of the $\text{P1-RT0-P0}^\text{DG}$ hierarchy of spaces, so we will refer to them
as compound $\text{P1-RT0-P0}^\text{DG}$ elements. Like the non-compound spaces described
by \citet{christiansen2008}, the expansion coefficients for $\mathbf{V}_0$
correspond to mesh vertices, for $\mathbf{V}_1$ to edges, and for $\mathbf{V}_2$
to cells. Thus, this hierarchy provides a finite element analogue of
the polygonal C-grid if we choose to represent velocity
in $\mathbf{V}_1$ and the mass variable in $\mathbf{V}_2$.

The construction of basis elements for the dual spaces proceeds in a very similar
way, except that the basis function for $\mathbf{V^1}$ is given by
$\mathbf{k} \times$ the solution of \refeqn{harmonicV1a} and \refeqn{harmonicV1c}.
This gives rise to a compound $\text{P1-N0-P0}^\text{DG}$ hierarchy of spaces,
where N0 refers to the lowest order two-dimensional N\'{e}d\'{e}lec elements.

\begin{figure}
\begin{centering}
\includegraphics[width=100mm]{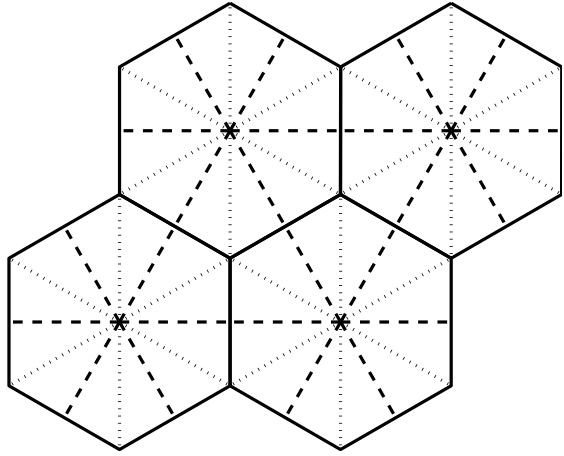}
\caption{
Example of part of a hexagonal primal mesh (solid lines) with its triangular
dual mesh (dashed lines) and the supermesh of triangular subelements (all
lines) used to construct the compound elements.}
\label{fig_supermesh}
\end{centering}
\end{figure}

An important detail concerns the number of subelements needed.
It may appear natural to subdivide an $n$-gon cell into $n$~triangular
subelements. However, it will be necessary to calculate integrals on the
overlap between primal and dual elements (section~\ref{massmatrices}). In order to be able to
do this when the domain is a curved surface approximated by plane triangular subelement
facets, both the primal and dual compound element meshes must be built
from triangular subelements of the the same supermesh. To achieve this
we divide $n$-gon cells (whether primal or dual) into $2n$~subelements
(Fig.~\ref{fig_supermesh}).

It is convenient to normalize the basis functions as follows:
\begin{eqnarray}
\alpha_i \in \mathbf{V}_2 :     & \ \ \ \ \ &
\int_{\mathrm{cell}\ i'} \alpha_i \, dA = \delta_{i \, i'} ; \\
\mathbf{v}_e \in \mathbf{V}_1 : & \ \ \ \ \ &
\int_{\mathrm{edge}\ e'} \mathbf{v}_e \cdot \mathbf{n} \, dl = \delta_{e \, e'} ; \\
\gamma_j \in \mathbf{V}_0 :     & \ \ \ \ \ &
\left. \gamma_j \right|_{\mathrm{vertex}\ j'} = \delta_{j \, j'} ; \\
\beta_j \in \mathbf{V}^2 :      & \ \ \ \ \ &
\int_{\mathrm{dual cell}\ j'} \beta_j \, dA = \delta_{j \, j'} ; \\
\mathbf{w}_e \in \mathbf{V}^1 : & \ \ \ \ \ &
\int_{\mathrm{dual edge}\ e'} \mathbf{w}_e \cdot \mathbf{m} \, dl = \delta_{e \, e'} ; \\
\chi_i \in \mathbf{V}^0 :       & \ \ \ \ \ &
\left. \chi_i \right|_{\mathrm{dual\ vertex}\ i'} = \delta_{i \, i'} .
\end{eqnarray}
Here $\mathbf{n}$ is the unit normal vector to primal edge $e$ and $\mathbf{m}$
is the unit tangent vector to dual edge $e$, 
with $\mathbf{m}$ and $\mathbf{n}$ pointing in the same sense
(i.e.\ $\mathbf{n} \cdot \mathbf{m} > 0$,
though they need not be parallel if the dual edges are not orthogonal
to the primal edges), as in \citet{thuburn2012}.
The normalization is chosen so that degrees of freedom for fields in
$\mathbf{V}_2$ and $\mathbf{V}^2$ correspond to area integrals of scalars
over primal cells and dual cells, respectively, degrees of freedom
for a field in $\mathbf{V}_1$ correspond to normal fluxes integrated along
primal edges, degrees of freedom in $\mathbf{V}^1$ correspond to circulations
integrated along dual edges, and degrees of freedom for fields in
$\mathbf{V}_0$ and $\mathbf{V}^0$ correspond to nodal values of scalars at primal
vertices and dual vertices respectively. Again, this corresponds closely
to the framework of \citet{thuburn2012}.

\begin{figure}
\begin{centering}

\includegraphics[width=50mm]{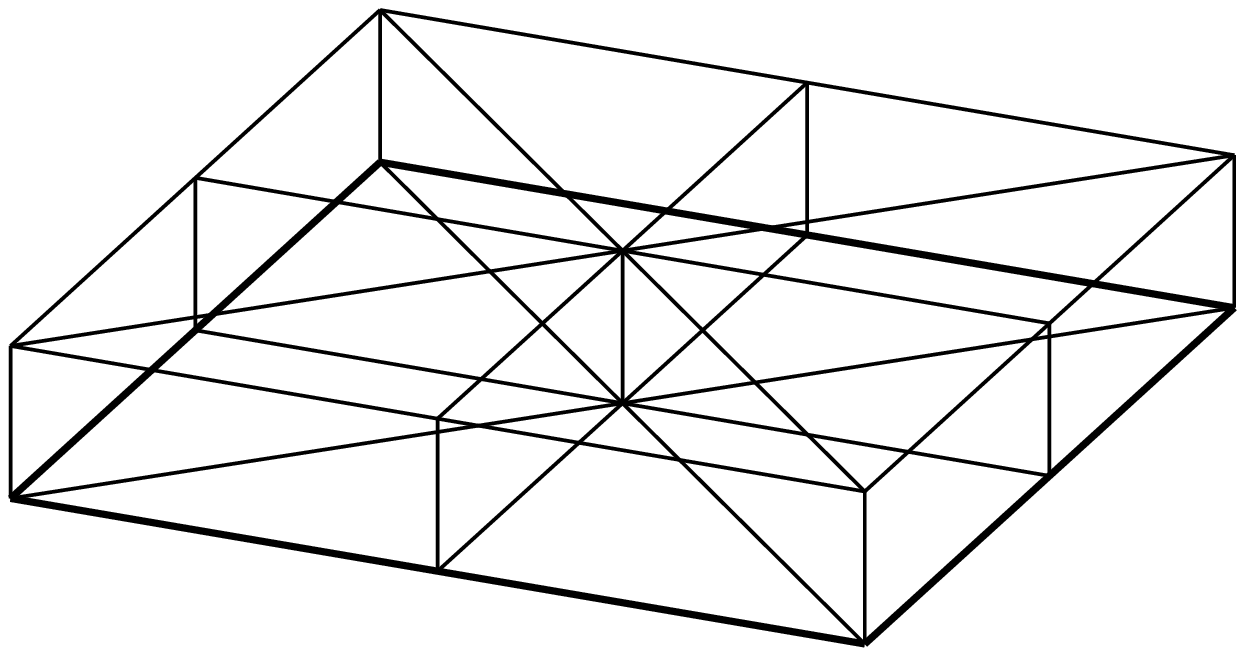}
~~~~~~~~~~
\includegraphics[width=50mm]{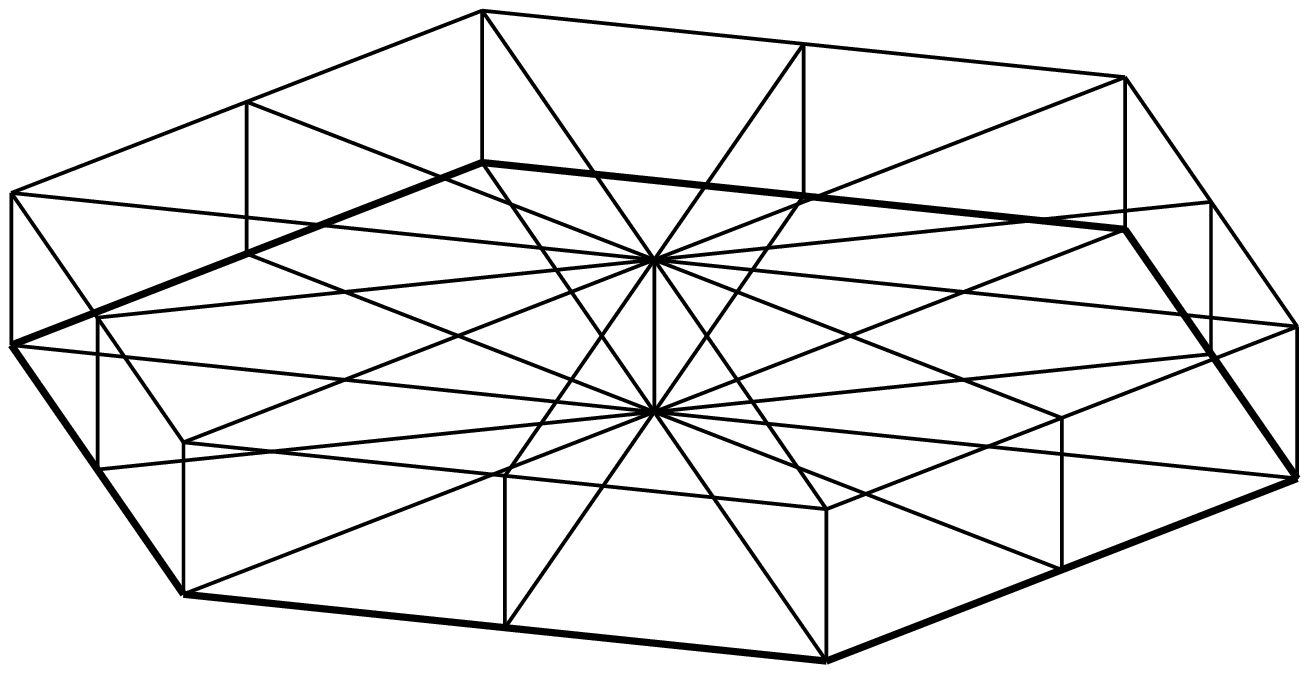}

\includegraphics[width=60mm]{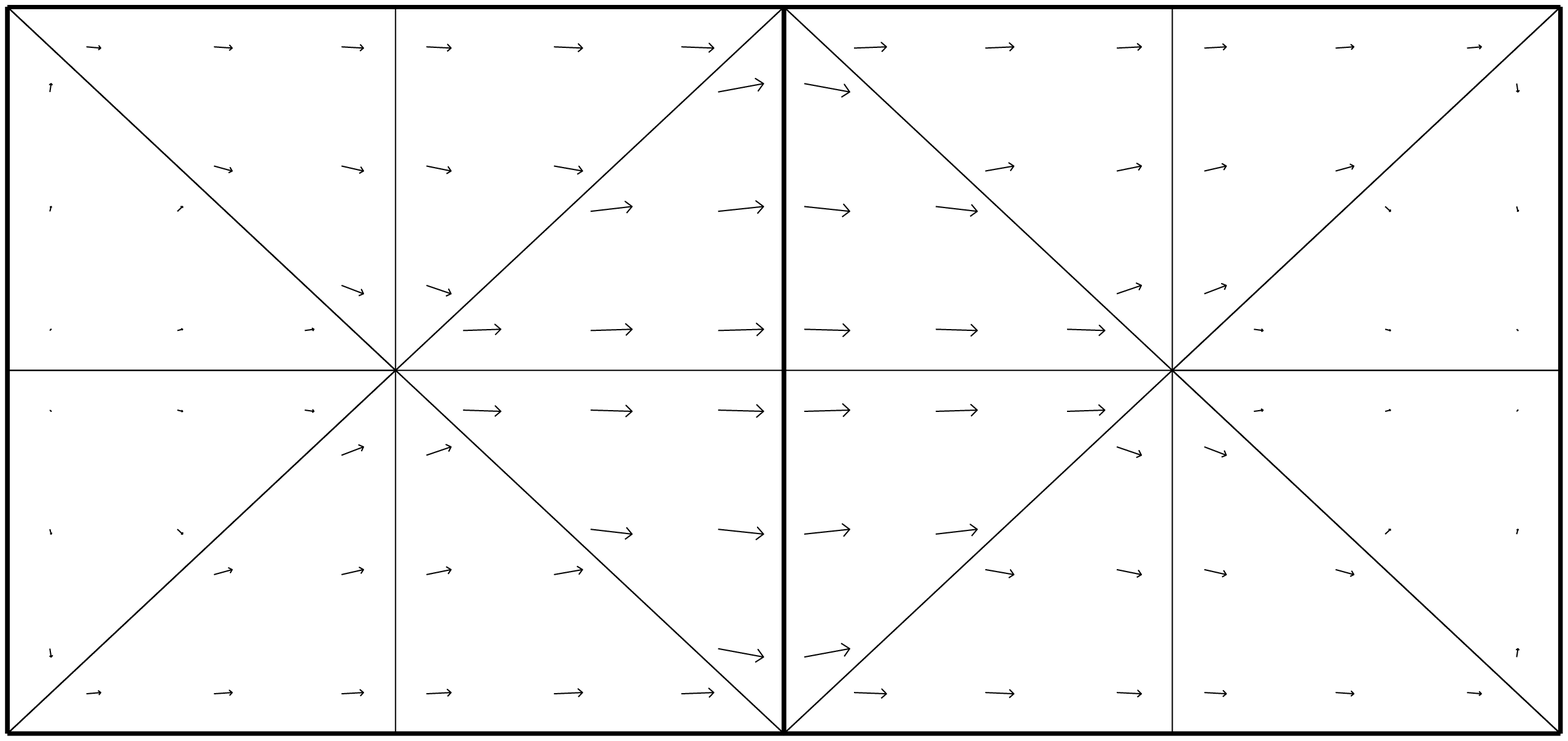}
\includegraphics[width=70mm]{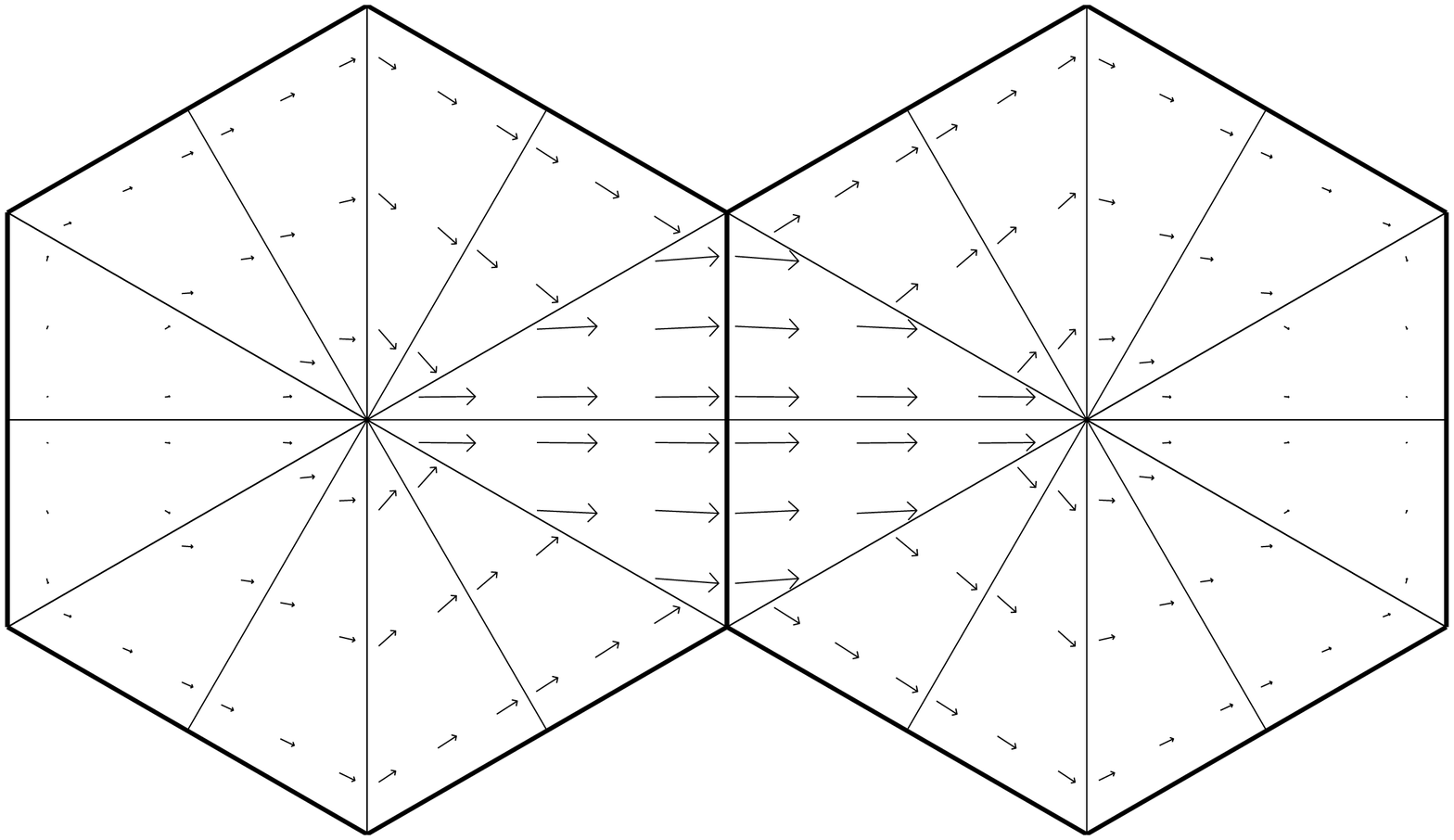}

\includegraphics[width=60mm]{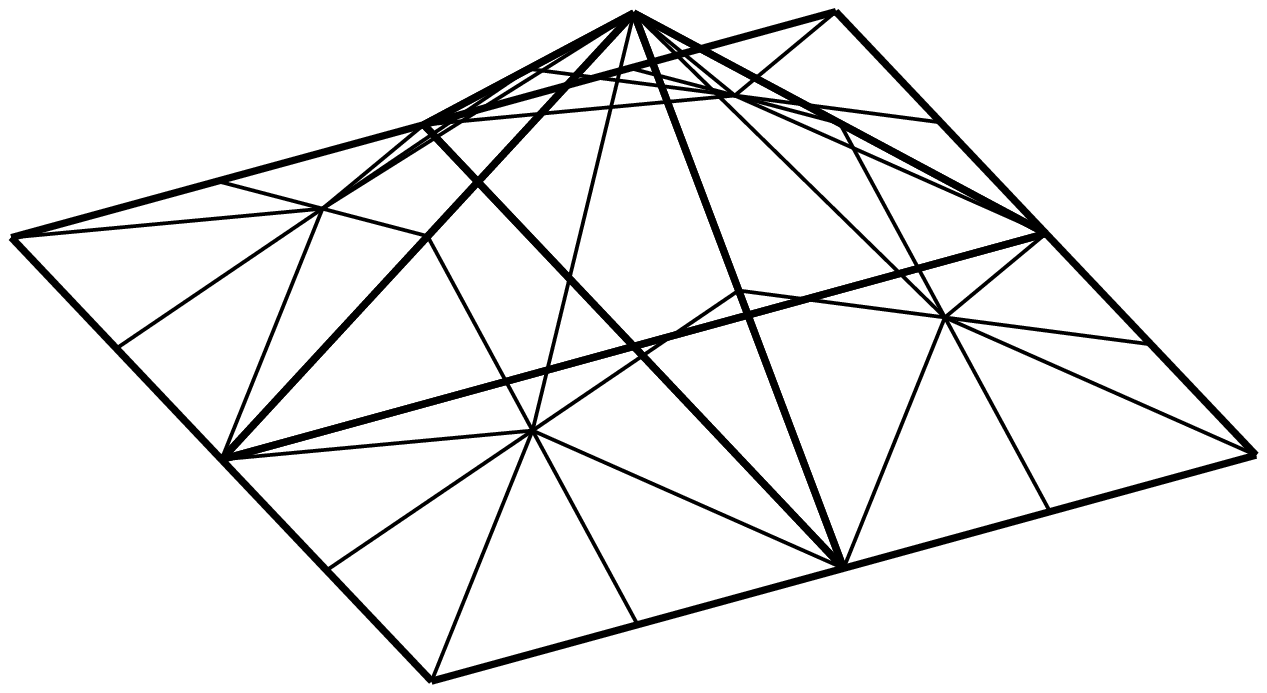}
~~~~
\includegraphics[width=60mm]{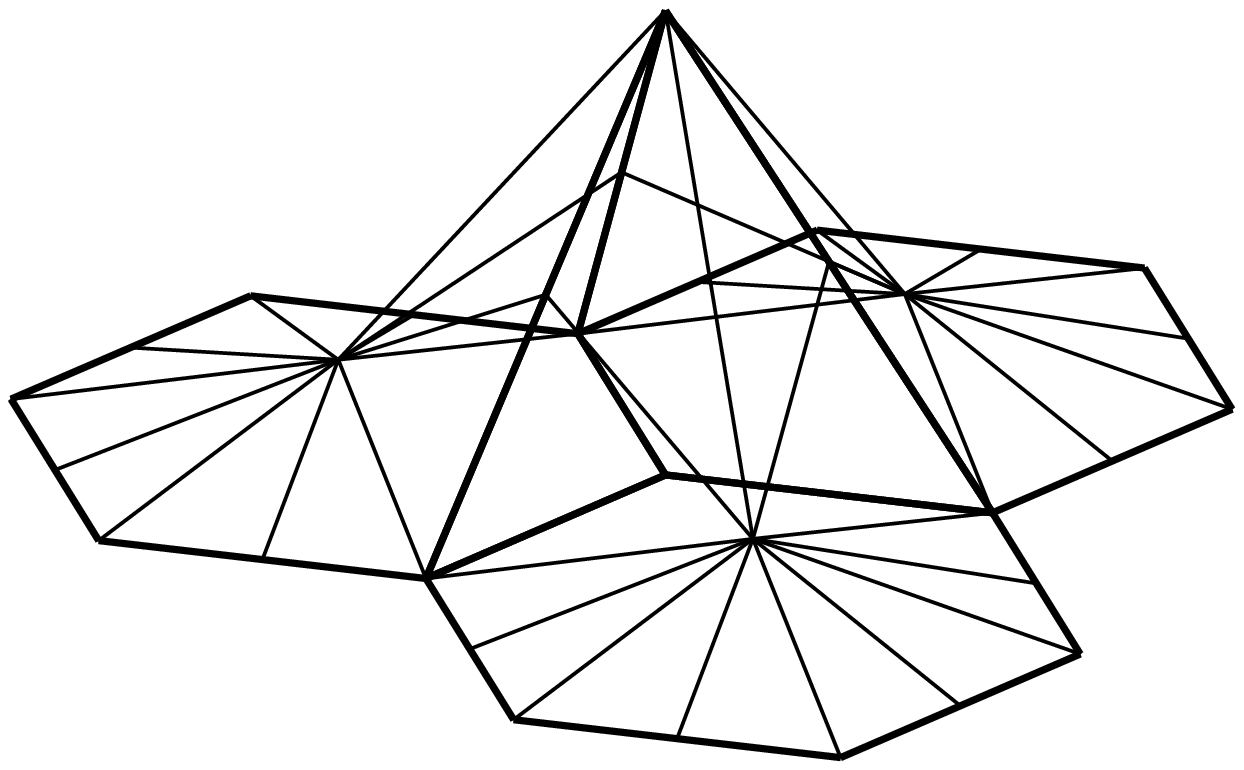}

\caption{
Typical compound basis elements of the function spaces on
a square mesh (left) and a hexagonal mesh (right).
structures.
Top row $\alpha_i \in \mathbf{V}_2$;
middle row $\mathbf{v}_e \in \mathbf{V}_1$;
bottom row $\gamma_j \in \mathbf{V}_0$.
In the middle row, at subelement edges the normal components
of the basis vectors are continuous.
}
\label{fig_elements}
\end{centering}
\end{figure}

\citet{melvin2014} have analyzed the wave dispersion properties for finite
element discretizations of the linear shallow water equations using these
compound elements. That paper gives explicit expressions for the 
$\mathbf{V}_1$ and $\mathbf{V}_2$ compound element basis functions
for the cases of a square mesh and a regular hexagonal mesh on a plane.
For more general meshes it is straightforward and convenient to construct
the compound element basis functions numerically.
Figure~\ref{fig_elements} shows typical basis elements for the three
primal mesh function spaces for quadrilateral and hexagonal cells.

The fields used in the computation are represented as expansions in terms of
these basis elements. For example,
\begin{eqnarray}
\phi = \sum_i \phi_i \alpha_i        & & \ \ \ \ \ \in \mathbf{V}_2, \label{eqn_exp_v2}\\
\mathbf{u} = \sum_e u_e \mathbf{v}_e & & \ \ \ \ \ \in \mathbf{V}_1 \label{eqn_exp_v1}
\end{eqnarray}
for the prognostic geopotential and velocity fields,
and
\begin{equation}
\xi = \sum_j \xi_j \gamma_j \ \ \ \ \ \in \mathbf{V}_0  \label{eqn_exp_v0}
\end{equation}
for the relative vorticity field. Here the sums are global sums over
all basis elements in the relevant spaces.
In some cases it will be useful to introduce dual space representations
of fields; these will be indicated by a hat symbol where necessary to
distinguish them from the corresponding primal space representations.
For example,
\begin{eqnarray}
\hat{\phi} = \sum_i \hat{\phi}_i \chi_i          & & \ \ \ \ \ \in \mathbf{V}^0, \label{eqn_exp_V0}\\
\hat{\mathbf{u}} = \sum_e \hat{u}_e \mathbf{w}_e & & \ \ \ \ \ \in \mathbf{V}^1, \label{eqn_exp_V1}\\
\hat{\xi} = \sum_j \hat{\xi}_j \beta_j           & & \ \ \ \ \ \in \mathbf{V}^2.  \label{eqn_exp_V2}
\end{eqnarray}
The fields $\phi$ and $\hat{\phi}$ have the same number of degrees of freedom,
and it is possible construct a well-conditioned and reversible map
between them by demanding that they agree when integrated
against any test function in the primal space $\mathbf{V}_2$.
Similarly, the fields $\mathbf{u}$ and $\hat{\mathbf{u}}$ have the same number
of degrees of freedom, and it is possible construct a well-conditioned and
reversible map between them by demanding that they agree when integrated
against any test function in the primal space $\mathbf{V}_1$.
It will also be useful to introduce spatially averaged versions of some fields.
For example,
\begin{eqnarray}
\widetilde{\phi} = \sum_j \widetilde{\phi}_j \gamma_j & & \ \ \ \ \ \in \mathbf{V}_0, \label{eqn_exp_v0_bar}\\
\overline{\phi} = \sum_j \overline{\phi}_j \beta_j    & & \ \ \ \ \ \in \mathbf{V}^2. \label{eqn_exp_V2_tilde}
\end{eqnarray}
Here, $\widetilde{\phi}$ and $\overline{\phi}$ have the same number of degrees of freedom,
and it is possible construct a well-conditioned and reversible map
between them by demanding that they agree when integrated
against any test function in the primal space $\mathbf{V}_0$. 
$\widetilde{\phi}$ or $\overline{\phi}$ can be obtained from $\phi$ by demanding
that they agree when integrated against any test function in $\mathbf{V}_0$; in effect this
provides an averaging operation from $\mathbf{V}_2$ to $\mathbf{V}_0$ or $\mathbf{V}^2$.
(However, we should not expect to be able to obtain $\phi$ from
$\widetilde{\phi}$ or $\overline{\phi}$,
as this would require an un-averaging operation, which will be ill-conditioned if
it exists at all.)

It will be convenient to be able to refer to the vector of degrees of
freedom for any field. To do this, we will use the same letter
(with hat, tilde or bar if needed) but in upper case. Thus, for example,
$\Phi$ will be the vector of values $(\phi_1, \phi_2, \dots )^T$, $\hat{U}$
will be the vector of values $( \hat{u}_1 , \hat{u}_2 , \ldots )^T$, etc.

\section{Finite element scheme}
\label{sec_FE_scheme}

Finite element schemes solve the governing equations by approximating the
solution in the chosen function spaces, written as expansions in
terms of basis functions (e.g.\ \refeqn{eqn_exp_v2}, \refeqn{eqn_exp_v1}),
and demanding that the equations be satisfied in weak form, that is,
when multiplied by any test function in the appropriate space and integrated
over the domain. In this approach \refeqn{ctsmass} becomes
\begin{equation}
\int \alpha_i \left( \phi_t + \nabla \cdot \mathbf{f} \right) \, dA
  =  0 \ \ \ \ \forall \alpha_i \in \mathbf{V}_2 ,
\end{equation}
or, regarding the integral as an inner product for which we introduce
angle backet notation,
\begin{equation}
\label{weakmass}
\langle \alpha_i ,  \phi_t \rangle
+
\langle \alpha_i , \nabla \cdot  \mathbf{f}  \rangle
= 0  \ \ \ \ \forall \alpha_i \in \mathbf{V}_2 .
\end{equation}
Similarly, \refeqn{ctsvel} becomes
\begin{equation}
\int \mathbf{v}_e \cdot
\left\{ \mathbf{u}_t + \mathbf{q}^{\perp} + \nabla (\phi_{\mathrm{T}} + k) \right\} \, dA
=  0 \ \ \ \ \forall \mathbf{v}_e \in \mathbf{V}_1 ,
\end{equation}
or
\begin{equation}
\label{weakvel}
\langle \mathbf{v}_e , \mathbf{u}_t \rangle
+
\langle \mathbf{v}_e , \mathbf{q}^{\perp} \rangle
+
\langle \mathbf{v}_e , \nabla (\phi_{\mathrm{T}} + k) \rangle = 0 .
\end{equation}
(The construction of the nonlinear terms $\mathbf{f}$, $\mathbf{q}$ and $k$ is
discussed in section~\ref{nonlinearswe} below.)
The method generally leads to a system of algebraic equations
for the unknown coefficients in the expansion of the solution.

The following subsections show how the mimetic finite element method can be
re-expressed in terms of certain matrix operators acting on the
coefficient vectors $\Phi$, $U$, etc. The notation is chosen to highlight
the similarity to the finite volume scheme of \citet{thuburn2014}.

\subsection{Matrix representation of derivatives -- strong derivatives}
\label{deriv1}

The velocity basis elements are constructed and normalized so as to have
constant divergence over the cell upwind of the edge where the degree of
freedom resides, with area integral equal to $1$, and constant
divergence over the cell downwind of this edge, with area integral
equal to $-1$, with zero velocity and hence zero divergence 
in all other cells. Thus
\begin{equation}
\nabla \cdot \mathbf{v}_e
=
\sum_i n_{e \, i} \alpha_i \ \ \ \ \ \in \mathbf{V}_2,
\end{equation}
where $n_{e \, i}$ is equal to $1$ when the normal at edge $e$
points out of cell $i$, equal to $-1$ when the normal at edge $e$ points
into cell $i$, and is zero otherwise. We will write $D_2$ for the
matrix whose transpose has components $n_{e \, i}$.
$D_2$ is called an incidence matrix because
it describes some aspects of the grid topology.
Hence, the divergence $\delta$ of an arbitrary velocity field $\mathbf{u}$ is
\begin{equation}
\sum_i \delta_i \alpha_i = \delta = \nabla \cdot \mathbf{u}
= \sum_e u_e \nabla \cdot \mathbf{v}_e
= \sum_{e \, i} u_e n_{e \, i} \alpha_i .   \label{eqn_div}
\end{equation}
Equating coefficients of $\alpha_i$ gives
\begin{equation}
\delta_i = \sum_e n_{e \, i} u_e,
\end{equation}
or, in matrix-vector notation
\begin{equation}
\Delta = D_2 U .
\end{equation}
Note we could have demanded that \refeqn{eqn_div} should hold when integrated
against any test function in $\mathbf{V}_2$, to obtain the same result.
However, this would obscure the fact that \refeqn{eqn_div} actually holds
at every point in the domain (except on cell edges where all terms are discontinuous),
not just when integrated against a test function.
In this sense, $\nabla \cdot \, : \mathbf{V}_1 \rightarrow \mathbf{V}_2$
is a \textit{strong} derivative operator.

Similarly, the basis elements in $\mathbf{V}_0$ are constructed
so that
\begin{equation}
\nabla^\perp \gamma_j
=
\sum_e - t_{e \, j} \mathbf{v}_e ,
\end{equation}
where $t_{e \, j}$ is defined to equal $1$ if edge $e$ is incident on
vertex $j$ and the unit
tangent vector $\mathbf{t}$ at edge $e$ points towards
vertex $j$, $-1$ if it points away from vertex $j$, and zero
otherwise. The unit normal and unit tangent at any edge are related by
$\mathbf{t} = \mathbf{k} \times \mathbf{n}$.
Hence, a stream function $\psi$ is related to the corresponding rotational
velocity field $\mathbf{u}$ by
\begin{equation}
\sum_e u_e \mathbf{v}_e
=
\mathbf{u}
=
\nabla^\perp \psi
=
\sum_j \psi_j \nabla^\perp \gamma_j
=
\sum_{j \, e} -\psi_j t_{e \, j} \mathbf{v}_e .    \label{eqn_gradperp}
\end{equation}
Equating coefficients of $\mathbf{v}_e$ and defining $D_1$ to be
the matrix whose entries are $t_{e \, j}$ gives the
matrix-vector form
\begin{equation}
U = - D_1 \Psi .
\end{equation}
Equation \refeqn{eqn_gradperp} holds pointwise (again with the exception of discontinuities),
so $\nabla^\perp : \mathbf{V}_0 \rightarrow \mathbf{V}_1$ is a strong derivative operator.

The matrices $D_1$ and $D_2$ are exactly the same as in the finite
volume framework of  \citet{thuburn2012}. In particular, they have the
property that
\begin{equation}
D_2 D_1 \equiv 0,
\end{equation}
giving a discrete analogue of the continuous property
$\nabla \cdot \nabla^\perp \equiv 0$.

Analogous relations hold on the dual spaces.
\begin{equation}
\nabla \chi_i
=
- \sum_e n_{e \, i} \mathbf{w}_e
\end{equation}
implies that the discrete analogue of
\begin{equation}
\hat{\mathbf{u}} = \nabla \hat{p}
\end{equation}
is
\begin{equation}
\hat{U} = \overline{D}_1 \hat{P} ,
\end{equation}
where $\overline{D}_1 = -D_2^T$.
Similarly
\begin{equation}
\mathbf{k} \cdot \nabla \times \mathbf{w}_e
=
\sum_j t_{e \, j} \beta_j
\end{equation}
implies that the discrete analogue of
\begin{equation}
\hat{\xi} = \mathbf{k} \cdot \nabla \times \hat{\mathbf{u}}
\end{equation}
is
\begin{equation}
\hat{\Xi} = \overline{D}_2 \hat{U},
\end{equation}
where $\overline{D}_2 = D_1^T$.
Again, these are strong derivative operators.

The matrices $\overline{D}_1$ and $\overline{D}_2$ have the property
\begin{equation}
\overline{D}_2 \overline{D}_1 \equiv 0,
\end{equation}
giving a discrete analogue in the dual space
of the continuous relation $\mathbf{k} \cdot \nabla \times \nabla \equiv 0$.

\subsection{Mass matrices and other operators}
\label{massmatrices}

Define the following mass matrices for the primal function spaces:
\begin{eqnarray}
\label{massL}
L_{i \, i'} = \langle \alpha_i , \alpha_{i'} \rangle
=
\int \alpha_i \alpha_{i'} \, dA ,
& \ \ \ \ \ & (\mathbf{V}_2 \rightarrow \mathbf{V}_2), \\
\label{massM}
M_{e \, e'} = \langle \mathbf{v}_e , \mathbf{v}_{e'} \rangle
=
\int \mathbf{v}_e \cdot \mathbf{v}_{e'} \, dA ,
& \ \ \ \ \ & (\mathbf{V}_1 \rightarrow \mathbf{V}_1), \\
\label{massN}
N_{j \, j'} = \langle \gamma_j , \gamma_{j'} \rangle
=
\int \gamma_j \gamma_{j'} \, dA ,
& \ \ \ \ \ & (\mathbf{V}_0 \rightarrow \mathbf{V}_0).
\end{eqnarray}
The expressions in parentheses indicate that $L$ maps
$\mathbf{V}_2$ to itself, etc.
(Analogous mass matrices may be defined for the dual spaces;
however, they will not be needed here.)

The following matrices are also needed.
\begin{eqnarray}
\label{operR}
R_{j \, i} = \langle \gamma_j , \alpha_{i} \rangle ,
& \ \ \ \ \ & (\mathbf{V}_2 \rightarrow \mathbf{V}_0), \\
\label{operW}
W_{e \, e'} = - \langle \mathbf{v}_e , \mathbf{v}_{e'}^\perp \rangle
= - W_{e' \, e} ,
& \ \ \ \ \ & (\mathbf{V}_1 \rightarrow \mathbf{V}_1), \\
H_{e \, e'} = \langle \mathbf{v}_e , \mathbf{w}_{e'} \rangle ,
& \ \ \ \ \ & (\mathbf{V}^1 \rightarrow \mathbf{V}_1), \\
J_{j \, j'} = \langle \gamma_j , \beta_{j'} \rangle ,
& \ \ \ \ \ & (\mathbf{V}^2 \rightarrow \mathbf{V}_0).
\end{eqnarray}
For completeness we may also define
\begin{equation}
I_{i \, i'} = \langle \alpha_i , \chi_{i'} \rangle ,
\ \ \ \ \ (\mathbf{V}^0 \rightarrow \mathbf{V}_2),
\end{equation}
though we will not need to employ this matrix in the shallow water scheme.

One further operator will be needed to construct the kinetic energy per unit mass.
It is
\begin{equation}
T_{i \, e \, e'} = \int_{\mathrm{cell}\ i} \mathbf{v}_e \cdot \mathbf{v}_{e'} \, dA
= A_i \langle \alpha_i , \mathbf{v}_e \cdot \mathbf{v}_{e'} \rangle
\ \ \ \ \ (\mathbf{V}_1 \otimes \mathbf{V}_1 \rightarrow \mathbf{V}_2).
\end{equation}
where $A_i = (L_{i \, i})^{-1}$ is the area of primal cell~$i$.

All of these matrices can be precomputed, so that no quadrature
needs to be done at run time. Moreover, they are all sparse, so they can be
efficiently stored as lists of stencils and coefficients.

Let $U^\perp$ be the coefficients of the expansion of the projection
of $\mathbf{u}^\perp$ into $\mathbf{V}_1$:
\begin{equation}
\label{uperp}
\langle \mathbf{v}_e , \sum_{e'} U^\perp_{e'} \mathbf{v}_{e'} \rangle
=
\langle \mathbf{v}_e , \sum_{e'} U_{e'} \mathbf{v}^\perp_{e'} \rangle
=
\langle \mathbf{v}_e , \mathbf{u}^\perp \rangle
~~~~~~~~\forall \mathbf{v}_e \in \mathbf{V}_1 .
\end{equation}
Using \refeqn{massM} and \refeqn{operW} gives the discrete version of the $\perp$
operator:
\begin{equation}
\label{perp}
M U^\perp = - W U .
\end{equation}

Demanding agreement between \refeqn{eqn_exp_v2} and \refeqn{eqn_exp_V0}
when integrated against any test function in $\mathbf{V}_2$ leads to
\begin{equation}
L \Phi = I \hat{\Phi} .
\label{eqn_hodge02}
\end{equation}
Similarly, demanding agreement between \refeqn{eqn_exp_v1} and
\refeqn{eqn_exp_V1} when integrated against any test function
in $\mathbf{V}_1$ gives
\begin{equation}
M U = H \hat{U} ,
\label{eqn_hodge11}
\end{equation}
while demanding agreement between  \refeqn{eqn_exp_v0} and
\refeqn{eqn_exp_V2} when integrated against any test function
in $\mathbf{V}_0$ gives
\begin{equation}
N \Xi = J \hat{\Xi} .
\label{eqn_hodge20}
\end{equation}
The relations \refeqn{eqn_hodge02}, \refeqn{eqn_hodge11}, \refeqn{eqn_hodge20}
provide invertible maps between the primal and dual
function spaces. 
Thus, they are examples of discrete Hodge star operators
\citep[e.g.][]{hiptmair2001}.
They may be contrasted with the analogous relations employed by
\citet{thuburn2012} and \citet{thuburn2014} for the finite volume case,
which do not involve mass matrices.

Demanding agreement between \refeqn{eqn_exp_v2}, \refeqn{eqn_exp_v0_bar},
and \refeqn{eqn_exp_V2_tilde} when integrated against any test function in
$\mathbf{V}_0$ leads to
\begin{equation}
N \widetilde{\Phi} = J \overline{\Phi} = R \Phi .
\label{eqn_avephi}
\end{equation}
This is the matrix representation of the averaging operator discussed in
section~\ref{sec_compound}.

\subsection{Matrix representation of derivatives -- weak derivatives}
\label{deriv2}

A field in $\mathbf{V}_2$ is discontinuous, so its gradient in
$\mathbf{V}_1$ can only be defined in a \textit{weak} sense,
by integrating against all test functions in $\mathbf{V}_1$.
For example,
\begin{equation}
\mathbf{g} = \nabla \phi
\end{equation}
must be approximated as
\begin{equation}
\langle \mathbf{v}_e , \mathbf{g} \rangle
=
\langle \mathbf{v}_e , \nabla \phi \rangle
\ \ \forall \mathbf{v}_e \in \mathbf{V}_1,
\end{equation}
where $\phi \in \mathbf{V}_2$, $\mathbf{g} \in \mathbf{V}_1$.
Expanding both $\phi$ and $\mathbf{g}$ in terms of basis elements
and integrating by parts then leads to the matrix form
\begin{equation}
M G = \overline{D}_1 L \Phi.
\end{equation}

Similarly, the curl of a vector field in $\mathbf{V}_1$ must be defined by
integration against all test functions in $\mathbf{V}_0$.
For example, the discrete analogue of
\begin{equation}
\xi = \mathbf{k} \cdot \nabla \times \mathbf{u},
\end{equation}
after expanding in basis functions and integrating by parts,
is
\begin{equation}
N \Xi = \overline{D}_2 M U. \label{xi}
\end{equation}

Combining these two results, the discrete analogue of
\begin{equation}
z = \mathbf{k} \cdot \nabla \times \nabla \phi
\end{equation}
is
\begin{equation}
N Z = \overline{D}_2 \overline{D}_1 L \Phi ,
\end{equation}
which is identically zero.

\begin{figure}
\begin{centering}
\makebox{\resizebox{150mm}{!}{{\includegraphics{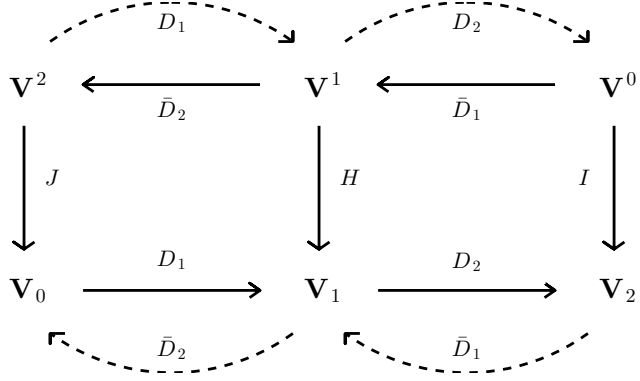}}}}
\caption{
Schematic showing the function spaces used in the scheme
and how the various matrices introduced above map between them.}
\label{opermap}
\end{centering}
\end{figure}

These derivative operators can be combined to obtain the
Laplacian of a scalar. For a scalar $\phi \in \mathbf{V}_2$,
the discrete Laplacian is $D_2 M^{-1} \overline{D}_1 L \Phi$.
For a scalar $\psi \in \mathbf{V}_0$, the discrete Laplacian
is $ - N^{-1} \overline{D}_2 M D_1 \Psi $.
The operators introduced above lead to a discrete version of the
Helmholtz decomposition, in which an arbitrary vector field is
decomposed into its divergent and rotational parts:
\begin{equation}
U = M^{-1} \overline{D}_1 L \Phi - D_1 \Psi .
\end{equation}
Figure~\ref{opermap} summarizes how the operators introduced here
map between the different function spaces.

\subsection{Some operator identities}
\label{identities}

The operators defined above satisfy some key relations that
underpin the mimetic properties of the scheme.
We have already seen that
\begin{equation}
{D}_2 {D}_1 \equiv 0
\end{equation}
and
\begin{equation}
\overline{D}_2 \overline{D}_1 \equiv 0,
\end{equation}
leading to discrete analogues of $\nabla \cdot \nabla^\perp \equiv 0$
and $\mathbf{k} \cdot \nabla \times \nabla \equiv 0$.

Next, note that the basis elements $\gamma_j$ give a partition of unity,
that is
\begin{equation}
\sum_j \gamma_j = 1
\end{equation}
at every point in the domain. Consequently
\begin{equation}
\label{sumR}
\sum_j R_{j \, i} = \langle 1 , \alpha_i \rangle = 1
\end{equation}
and
\begin{equation}
\sum_j J_{j \, j'} = \langle 1 , \beta_{j'} \rangle = 1.
\end{equation}

Now let
\begin{equation}
\psi = \sum_j \psi_j \gamma_j \ \ \in \mathbf{V}_0 .
\end{equation}
By considering the projection of $\nabla \psi$ into $\mathbf{V}_1$
\begin{equation}
\langle \mathbf{v}_e , (\nabla^\perp \psi)^\perp \rangle =
    - \langle \mathbf{v}_e , \nabla \psi \rangle
\end{equation}
and integrating by parts and using the matrices defined in
sections~\ref{deriv1} and~\ref{massmatrices}, we obtain
\begin{equation}
- \overline{D}_2 W = R D_2. \label{TRiSK}
\end{equation}
This identity is key to obtaining the steady geostrophic mode property
(section~\ref{geostrophic} below).
A rough interpretation is that averaging velocities to construct the
Coriolis terms ($W$) then taking their divergence ($\overline{D}_2$)
gives the same result as computing the velocity divergence ($D_2$)
followed by averaging to $\mathbf{V}_0$ ($R$). 
One consequence is $\overline{D}_2 W D_1 = - R D_2 D_1 \equiv 0$.

An identical formula to \refeqn{TRiSK}
relating $R$ and $W$ was obtained by \citet{thuburn2012} for the finite volume case.
The result was originally derived for the construction of the Coriolis
terms on orthogonal grids by \citet{thuburn2009}, and \citet{thuburn2012}
showed that it could be embedded in a more general framework applicable
to nonorthogonal grids. Moreover, \citet{thuburn2009} showed that, for
any given $R$ with the appropriate stencil (which we have here)
and satisfying \refeqn{sumR},
there is a unique antisymmetric $W$ satisfying \refeqn{TRiSK}, and gave an
explicit construction for $W$ in terms of $R$. Thus, although the context
and interpretation are slightly different here, we can, nevertheless, use
the Thuburn et al.\ construction in implementing the mixed finite-element
version of the $W$ operator!

Now consider the two representations of any vector field
$\mathbf{u} \in \mathbf{V}_1$, $\hat{\mathbf{u}} \in \mathbf{V}^1$
related by
\begin{equation}
\langle \mathbf{v} , \mathbf{u} \rangle
=
\langle \mathbf{v} , \hat{\mathbf{u}} \rangle \ \ \forall \mathbf{v} \in \mathbf{V}_1,
\end{equation}
so that
\begin{equation}
\label{uuhat}
M U = H \hat{U}.
\end{equation}
Since $\nabla^\perp \gamma_j \in \mathbf{V}_1$,
\begin{equation}
\langle \nabla^\perp \gamma_j , \mathbf{u} \rangle
=
\langle \nabla^\perp \gamma_j , \hat{\mathbf{u}} \rangle 
 \ \ \forall \gamma_j \in \mathbf{V}_0,
\end{equation}
and integrating by parts gives
\begin{equation}
- \langle \gamma_j , \mathbf{k} \cdot \nabla \times \mathbf{u} \rangle
=
- \langle \gamma_j , \mathbf{k} \cdot \nabla \times \hat{\mathbf{u}} \rangle
 \ \ \forall \gamma_j \in \mathbf{V}_0.
\end{equation}
Hence
\begin{equation}
- \overline{D}_2 M U = - J \overline{D}_2 \hat{U} .
\end{equation}
Finally, substituting from \refeqn{uuhat} and noting that 
$\hat{U}$ is arbitrary gives
\begin{equation}
\label{JDeqDH}
\overline{D}_2 H = J \overline{D}_2 .
\end{equation}
The interpretation of this identity is that, for a velocity field in
$\mathbf{V}^1$, taking the curl followed by mapping to the primal space
is equivalent to mapping the velocity
field to the primal space then taking its curl. One consequence is that
$\overline{D}_2 H \overline{D}_1 = J \overline{D}_2 \overline{D}_1 \equiv 0$.

Finally, let $\chi \in \mathbf{V}_2$ and $\hat{\chi} \in \mathbf{V}^0$ be
two discrete representations of a scalar field related by
\begin{equation}
\langle \alpha , \chi \rangle = \langle \alpha , \hat{\chi} \rangle
\ \ \forall \alpha \in \mathbf{V}_2 ,
\end{equation}
so that
\begin{equation}
\label{xxhat}
L X = I \hat{X} .
\end{equation}
Since $\nabla \cdot \mathbf{v}_e \in \mathbf{V}_2$ for any
$\mathbf{v}_e \in \mathbf{V}_1$, we have
\begin{eqnarray}
\langle \nabla \cdot \mathbf{v}_e , \chi \rangle
& = &
\langle \nabla \cdot \mathbf{v}_e , \hat{\chi} \rangle , \nonumber \\
\langle \mathbf{v}_e , \nabla \chi \rangle
& = &
\langle \mathbf{v}_e , \nabla \hat{\chi} \rangle , \nonumber \\
\overline{D}_1 L X & = & H \bar{D}_1 \hat{X},
\end{eqnarray}
or, using \refeqn{xxhat} and noting that $\hat{X}$ is arbitrary,
\begin{equation}
\label{DIeqHD}
\overline{D}_1 I = H \overline{D}_1 .
\end{equation}

Using these identities and the Hodge star operators, it can be seen
that taking a weak derivative in the primal space is equivalent to
applying a Hodge star to map to the dual space, taking a strong
derivative in the dual space, and applying another Hodge star to
map back to the primal space:
\begin{eqnarray}
M^{-1} \overline{D}_1 L & = & (M^{-1} H) \, \overline{D}_1 \, (I^{-1} L) ; \\
N^{-1} \overline{D}_2 M & = & (N^{-1} J) \, \overline{D}_2 \, (H^{-1} M) 
\end{eqnarray}
\citep{cotter2014}.
Thus, certain paths in Fig.~\ref{opermap} commute.
Weak derivative operators in the dual space can be defined by demanding
a similar equivalence with primal space strong derivatives; however,
the resulting formulas are less elegant and, in any case, will not
be needed.

\subsection{Linear shallow water equations}
\label{linearswe}

We first examine the spatial discretization of the linear shallow
water equations to illustrate how some key conservation and balance
properties arise.
The rotating shallow water equations \refeqn{ctsmass}, \refeqn{ctsvel}
when linearized about a resting
basic state with constant geopotential $\phi_0$ and with constant Coriolis
parameter $f$ become
\begin{eqnarray}
\phi_t + \nabla \cdot ( \phi_0 \mathbf{u} ) & = & 0 , \\
\mathbf{u}_t 
+ f \mathbf{u}^\perp + \nabla \phi & = & 0.
\end{eqnarray}
By writing these in weak form (analogous to \refeqn{weakmass} and \refeqn{weakvel}),
expanding $\phi$ and $\mathbf{u}$ in terms of basis functions,
and using the notation and operators defined above, we obtain
\begin{eqnarray}
\dot{\Phi} + \phi_0 D_2 U & = & 0 ,    \\
M \dot{U} - f W U + \overline{D}_1 L \Phi & = & 0 .
\end{eqnarray}

\subsubsection{Mass conservation}
\label{linearmass}

Mass conservation is trivially satisfied (for both the linear and nonlinear
equations) because the discrete divergence is a strong operator, so the
domain integral of the discrete divergence of any vector field vanishes.

\subsubsection{Energy conservation}
\label{energy}

For the linearized equations the total energy is given by
\begin{eqnarray}
E & = &
\frac{1}{2} \int \phi^2 + \phi_0 \mathbf{u} \cdot \mathbf{u} \, dA \nonumber \\
& = &
\frac{1}{2} \Phi^T L \Phi + \frac{1}{2} \phi_0 U^T M U .
\end{eqnarray}
Hence, the rate of change of total energy is
\begin{eqnarray}
\frac{d E}{d t} & = &
\Phi^T L \dot{\Phi} + \phi_0 U^T M \dot{U} \nonumber \\
& = &
-\phi_0 \Phi^T L D_2 U + \phi_0 U^T (f W U - \overline{D}_1 L \Phi) \nonumber \\
& = & 0 ,
\end{eqnarray}
where we have used the fact that $L$ and $M$ are symmetric, $W$ is antisymmetric,
and $D_2^T = - \overline{D}_1$.

\subsubsection{Steady geostrophic modes}
\label{geostrophic}

The linear shallow water equations support steady non-divergent
flows in geostrophic balance. A numerical method must respect this property
in order to be able to represent geostrophic balance. However, it is
non-trivial to achieve this property because several ingredients must
fall into place.

\begin{itemize}

\item
The geopotential $\phi$ must be steady. The steadiness of $\phi$
follows immediately from the assumption that $\nabla \cdot \mathbf{u} = 0$.

\item
The relative vorticity $\xi$ must be steady; neither the pressure
gradient nor the Coriolis term should generate vorticity.
First note that, from Property List~1, $U = - D_1 \Psi$
for some $\Psi$. Taking the curl of the momentum equation then gives
\begin{equation}
N \dot{\Xi} = \overline{D}_2 M \dot{U}
= \overline{D}_2 (- fW D_1 \Psi - \overline{D}_1 L \Phi).
\label{linvort}
\end{equation}
The pressure gradient term does not contribute because
$\overline{D}_2 \overline{D}_1 \equiv 0$, and the Coriolis
term does not contribute because $\overline{D}_2 W D_1 \equiv 0$.

\item
There must exist a geopotential $\phi$ that balances the Coriolis
term so that the divergence is steady. Taking the divergence of the momentum
equation gives the divergence tendency
\begin{equation}
\dot{\Delta} = D_2 M^{-1} (- f W D_1 \Psi - \overline{D}_1 \Phi) .
\end{equation}
If we define $\Phi = f L^{-1} R^T \Psi$ and use the transpose of
\refeqn{TRiSK} we find that $\dot{\Delta}$ does indeed vanish; thus the
required $\phi$ does exist.

\end{itemize}
Consequently, the scheme does support steady geostrophic modes for the
linearized equations.
(Note, it is not necessarily true that any given $\phi$
field can be balanced by some non-divergent velocity field.
On some meshes, particularly those with triangular primal cells,
there might not be enough velocity degrees of freedom
to balance all possible $\phi$ fields.)

\subsubsection{Linear PV equation}
\label{linPV}

A generalization of the steady geostrophic mode property is that the
scheme should have a suitable PV equation. In this section we consider the
linear case; the nonlinear case is dealt with in section~\ref{nonlinearPV}.

The mass field $\phi$ and the vorticity field $\xi$ live in different
spaces. To construct a suitable discrete PV we need an averaged mass
field $\widetilde{\phi}$
that lives in the same space as $\xi$. The linearized PV should be
independent of time. For this to hold, $\xi$ and $\widetilde{\phi}$
must see the same divergence field.

For a general (possibly divergent) velocity field $\mathbf{u}$, the vorticity
equation \refeqn{linvort} becomes
\begin{equation}
N \dot{\Xi} = f \overline{D}_2 W U .
\end{equation}
Define $\widetilde{\phi}$ using \refeqn{eqn_avephi}.
Then the evolution of $\widetilde{\phi}$ is given by
\begin{eqnarray}
N \dot{\widetilde{\Phi}} & = & R \dot{\Phi} \nonumber \\
& = & - \phi_0 R D_2 U \nonumber \\
& = & \phi_0 \overline{D}_2 W U ,
\end{eqnarray}
(using \refeqn{TRiSK}). Thus, $\Xi$ and $\widetilde{\Phi}$
see the same divergence $- \overline{D}_2 W U$; consequently the
linearized PV $\Xi / \phi_0 - f \widetilde{\Phi} / \phi_0^2$
is independent of time.

\subsection{Nonlinear shallow water equations}
\label{nonlinearswe}

The nonlinear rotating shallow water equations are
\refeqn{ctsmass} and \refeqn{ctsvel}. Writing these in weak form
\refeqn{weakmass} and \refeqn{weakvel}, and
letting $F$, $Q$, and $K$ be the vectors of coefficients for the discrete
representations of the mass flux $\mathbf{f} = \mathbf{u} \phi$,
the PV flux $\mathbf{q} = \mathbf{f} \pi$,
and the kinetic energy per unit mass $k = \mathbf{u} \cdot \mathbf{u} / 2$,
the nonlinear discretization becomes
\begin{eqnarray}
\label{phidot}
\dot{\Phi} + D_2 F & = & 0 , \\
\label{mudot}
M \dot{U} + M Q^\perp + \overline{D}_1 L ( \Phi_T + K ) & = & 0 .
\end{eqnarray}
The remaining issue is how to construct suitable values of 
the three nonlinear terms $K$, $F$, and $Q^\perp$.

\subsubsection{Constructing $K$}

The discretization of $k$ follows the standard finite element construction,
which is to project $\nabla k$ into $\mathbf{V}_1$.
It may easily be verified that this is equivalent to projecting
$k$ into $\mathbf{V}_2$ before taking the weak gradient.
Using the $T$ operator defined in section~\ref{massmatrices}, the
expansion coefficients $K$ of the projected $k$ are given by
\begin{equation}
k_i = \frac{1}{2} \sum_{e \, e'} T_{i \, e \, e'} u_e u_{e'} .
\end{equation}

\subsubsection{Constructing $F$}

Because the $\phi$ field is approximated as piecewise constant,
its degrees of freedom can be interpreted as primal cell integrals.
Similarly, the degrees of freedom of the $\mathbf{u}$ field
are the integrals of the normal velocity fluxes across primal cell
edges, and the $D_2$ operator looks exactly like a finite volume
divergence operator. Thus, it is straightforward to use a finite
volume advection scheme for advection of $\phi$.
The mass flux is constructed using a forward in time advection
scheme, identical to that used by \citet{thuburn2014},
using the fluxes $U$ and the mass field $\Phi$ as input.
We write this symbolically as
\begin{equation}
\label{adv1}
F = \mathrm{adv}_1 (U,\Phi) .
\end{equation}
The subscript $1$ indicates that this version of the advection scheme
operates on the primal mesh and works with densities or concentrations.

\subsubsection{Constructing $Q^\perp$}
\label{nonlinearPV}

So far we have not needed to use the dual mesh representation of any field.
However, in order to use the same finite volume advection scheme as
\citet{thuburn2014} to compute PV fluxes, we need a piecewise constant
representation of the PV field on dual cells, and a representation of the mass flux
field in terms of components normal to dual cell edges. These are naturally
given by the dual function spaces:
\begin{equation}
\pi = (f + \hat{\xi}) / \bar{\phi}~~~~~~~~\in \mathbf{V^2} ,
\end{equation}
and
\begin{equation}
\widehat{\mathbf{f}^\perp}~~~~~~~~\in \mathbf{V}^1 .
\end{equation}
Applying \refeqn{perp} followed by \refeqn{eqn_hodge11} to the mass flux 
gives
\begin{equation}
\label{hatfperp}
H \widehat{F^\perp} = M F^\perp = -W F .
\end{equation}

Now consider the evolution of the dual mass field $\bar{\phi}$.
\begin{equation}
J \dot{\bar{\Phi}} =
N \dot{\widetilde{\Phi}} =
\overline{D}_2 W F =
- \overline{D}_2 H \widehat{F^\perp} =
- J \overline{D}_2 \widehat{F^\perp} ,
\end{equation}
i.e.\ 
\begin{equation}
\label{phibar}
\dot{\bar{\Phi}} + \overline{D}_2 \widehat{F^\perp} = 0 .
\end{equation}
Since $\overline{D}_2$ acts exactly like a finite volume divergence operator
on the dual mesh, $\bar{\phi}$ behaves exactly as if it were evolving
according to a finite volume advection scheme.

Next, in order for PV to evolve in a way consistent with the
mass field $\bar{\phi}$, we construct PV fluxes in $\mathbf{V}^1$ using the
dual mesh finite volume advection scheme:
\begin{equation}
\label{adv2}
\widehat{Q^\perp} = \mathrm{adv}_2 (\widehat{F^\perp} , \Pi)
\end{equation}
The subscript~2 indicates that this version of the advection scheme
operates on the dual mesh and works with quantities analogous to
mixing ratios (such as PV $\pi$).
Finally, these dual mesh PV fluxes are mapped to the primal mesh
for use in the momentum equation:
\begin{equation}
H \widehat{Q^\perp} = M Q^\perp .
\end{equation}

It may be verified that the resulting vorticity equation for $\dot{\hat{\Xi}}$
is indeed analogous to \refeqn{phibar}, involving the potential vorticity flux
$\widehat{Q^\perp}$. Using \refeqn{eqn_hodge20}, \refeqn{xi}, \refeqn{mudot} and 
\refeqn{JDeqDH}, we have
\begin{equation}
\begin{array}{ccccccc}
J \dot{\hat{\Xi}} & = &   N \dot{\Xi}
                  & = &   \overline{D}_2 M \dot{U}
                  & = & - \overline{D}_2 M Q^\perp \nonumber \\
                  & = & - \overline{D}_2 H \widehat{Q^\perp}
                  & = & - J \overline{D}_2 \widehat{Q^\perp} . & &
\end{array}
\end{equation}
Hence,
\begin{equation}
\label{dualPV}
\dot{\hat{\Xi}} + \overline{D}_2 \widehat{Q^\perp} = 0,
\end{equation}
which is of the desired form. The similarity of \refeqn{dualPV} and \refeqn{phibar}
means that it is possible to construct PV fluxes from the dual mass fluxes
$\widehat{F^\perp}$ such that the evolution of the PV is consistent with
the evolution of $\bar{\Phi}$.

\subsection{Time integation scheme}

The same time integration scheme as in \citet{thuburn2014} is used.
\begin{eqnarray}
\label{timephi}
\Phi^{n+1} - \Phi^n + D_2 \widetilde{F} & = & 0, \\
\label{timeu}
M U^{n+1} - M U^n + H \widetilde{\widehat{Q^\perp}}
+ \overline{\overline{D}_1 L (\Phi + K)}^t & = & 0.
\end{eqnarray}
Here, $\overline{()}^t$ indicates the usual (possibly off-centred)
Crank-Nicolson approximation to the integral over one time interval:
\begin{equation}
\overline{\psi}^t = (\alpha \psi^{n+1} + \beta \psi^n) \Delta t
\end{equation}
(for any field $\psi$)
where $\alpha + \beta = 1$. All results presented below use
$\alpha = \beta = 0.5$.

$\widetilde{F}$ is an approximation to the time integral of the mass flux
across primal cell edges computed using the advection scheme. The velocity
field used for the advection is $\overline{U}^t$. We write this symbolically
as
\begin{equation}
\widetilde{F} = \mathrm{Adv}_1 (\overline{U}^t,\Phi^n) .
\end{equation}
(The notation $\mathrm{Adv}_1$, as distinct from $\mathrm{adv}_1$ in \refeqn{adv1},
indicates that here we are working with time integrals $\overline{U}^t$ and
$\widetilde{F}$.)

Finally, $\widetilde{\widehat{Q^\perp}}$ is an approximation to
the time integral of the PV flux across dual edges computed using the
advection scheme. Dual grid time integrated mass fluxes are calculated
from the primal grid time integrated mass fluxes as
\begin{equation}
\label{eq123}
H \widetilde{\widehat{F^\perp}} = -W \widetilde{F} .
\end{equation}
These are then used in the dual grid advection scheme to compute the
time integrated PV fluxes:
\begin{equation}
\widetilde{\widehat{Q^\perp}} = \mathrm{Adv}_2 (\widetilde{\widehat{F^\perp}},\Pi^n) .
\end{equation}

\subsection{Incremental iterative solver}

The system \refeqn{timephi}, \refeqn{timeu} is nonlinear in the unknowns
$\Phi^{n+1}$, $U^{n+1}$. It can be solved efficiently using an incremental
method; this may be viewed as a Newton method with an approximate Jacobian.
After $l$ iterations \refeqn{timephi} and \refeqn{timeu} will not be satisfied
exactly but will have some residuals $R_\Phi$, $R_U$ defined by:
\begin{eqnarray}
\label{resphi}
R_\Phi & = & \Phi^{(l)} - \Phi^n + D_2 \widetilde{F} , \\
\label{resu}
R_U & = & M U^{(l)} - M U^n + H \widetilde{\widehat{Q^\perp}}
+ \overline{\overline{D}_1 L (\Phi + K)}^t .
\end{eqnarray}
Here $\Phi^{(l)}$ and $U^{(l)}$ are the approximations after $l$ iterations
to $\Phi^{n+1}$ and $U^{n+1}$ and it is understood that these have been used in
evaluating $\widetilde{F}$, $\widetilde{\widehat{Q^\perp}}$,
and $\overline{\overline{D}_1 L (\Phi + K)}^t$.
We then seek updated values
\begin{equation}
\label{increment}
\Phi^{(l+1)} = \Phi^{(l)} + \Phi', \ \ \ \ \ U^{(l+1)} = U^{(l)} + U',
\end{equation}
that will reduce the residuals, where the increments $\Phi'$, $U'$ satisfy
\begin{equation}
\Phi' + \alpha \Delta t D_2 \phi^* U' = -R_\Phi,
\end{equation}
\begin{equation}
\label{uinc}
U' + \alpha \Delta t \mathcal{M}^{-1} \overline{D}_1 L \Phi' = - \mathcal{M}^{-1} R_U.
\end{equation}
Here, $\phi^*$ is a reference value of $\phi$; in the current implementation
it is given by $\phi^n$ interpolated to cell edges.
To avoid the appearance of the non-sparse matrix $M^{-1}$
in the Helmholtz problem below, a sparse approximation
$\mathcal{M}^{-1}$ has been introduced. The construction of
$\mathcal{M}^{-1}$ is briefly discussed in the Appendix.

Eliminating $U'$ leaves a Helmholtz problem for $\Phi'$:
\begin{equation}
\alpha^2 \Delta t^2 D_2 \phi^* \mathcal{M}^{-1} \overline{D}_1 L \Phi'
- \Phi' =
R_\Phi - \alpha \Delta t D_2 \phi^* \mathcal{M}^{-1} R_U .
\end{equation}
In the current implementation, the Helmholtz problem is solved using
a single sweep of a full multigrid algorithm. This gives sufficient
accuracy to avoid harming the convergence rate of the Newton
iteration.
Once $\Phi'$ is found, $U'$ is obtained by backsubstitution in \refeqn{uinc}.
Finally, \refeqn{increment} is used to obtain improved estimates for the unknowns.

Testing to date has given satisfactory results with 4~Newton iterations.
The algorithm requires the inversion of several of the linear operators
represented as matrices above. The appendix describes how this is done.

\section{Results}
\label{sec_results}

The same tests were applied to the finite element shallow water model as
were applied to the finite volume model of \citet{thuburn2014}. Only a subset
of results are shown here to emphasize the differences between the two models.
Other aspects are the following.
\begin{itemize}
\item
{\bf Stability.} All experimentation to date suggests the two models have the same
stability limit: with no temporal off-centring ($\alpha = \beta = 0.5$) the models
are stable for large gravity wave Courant numbers and advective Courant numbers
less than~1.

\item
{\bf Advection.} The same advection scheme is used in the two models to compute
mass, PV, and tracer fluxes on primal and dual meshes. In particular, the
models share the consistency between mass and PV, between mass and tracers,
and between primal mass and dual mass discussed by \citet{thuburn2014}.

\item
{\bf Balance.} The balance test discussed in section~6.8 of \citet{thuburn2014}
was repeated for the finite element shallow water model. The results on both
the hexagonal and cubed sphere meshes were very similar to those for the
finite volume model and the ENDGame semi-implicit semi-Lagrangian model
\citep{zerroukat2009},
implying that any spurious numerical generation of imbalance is extremely weak.

\item
{\bf Computational Rossby modes.} The experiment to test the ability of the scheme
on hexagonal meshes to handle grid-scale vorticity features was not repeated
here. However, given the general arguments in \citet{thuburn2014}
\citep[see also][]{weller2012}, and the similarities between the numerics of the
finite volume and finite element models, the results are expected to be very
similar for the finite element model.

\end{itemize}
For the remaining tests discussed below, the same mesh resolutions and time steps
were used as in \citet{thuburn2014}.

\subsection{Convergence of the Laplacian}

The discrete Laplacian defined in section~\ref{deriv2} was applied to the $\mathbf{V}_2$
representation of the field $\cos \varphi \sin \lambda$ on the unit sphere,
where $\varphi$ is latitude and $\lambda$ is longitude, and the $L_\infty$ and $L_2$
errors computed on different resolution meshes. The results are shown in table~\ref{tab:laplacian}.

\begin{table}
\centering
\caption{Convergence of the scalar Laplacian on hexagonal and cubed sphere grids.}
\label{tab:laplacian}
\begin{tabular}{|r|l|l|r|l|l|}
\hline
Hex    &                 &            & Cube   &                &           \\
Ncells & $L_\infty$ err  &  $L_2$ err & Ncells & $L_\infty$ err & $L_2$ err \\
\hline
    42 & 0.14            & 0.074      &     54 & 0.12           & 0.064     \\
   162 & 0.033           & 0.019      &    216 & 0.030          & 0.016     \\
   642 & 0.0090          & 0.0049     &    864 & 0.0077         & 0.0043    \\
  2562 & 0.0026          & 0.0012     &   3456 & 0.0038         & 0.0012    \\
 10242 & 0.00082         & 0.00031    &  13824 & 0.0022         & 0.00037   \\
 40962 & 0.00036         & 0.000081   &  55296 & 0.0012         & 0.00012   \\
163842 & 0.00018         & 0.000022   & 221184 & 0.00062        & 0.000039  \\
\hline
\end{tabular}
\end{table}

On both the hexagonal and cubed sphere meshes the $L_\infty$ errors converge
at first order. On the hexagonal mesh the $L_2$ errors converge at close to
second order, while on the cubed sphere mesh the convergence rate is between
first and second order. For the cubed sphere mesh the convergence of the discrete
scalar Laplacian is significantly better than for the finite volume scheme
of \citet{thuburn2014} (their table~4).

\subsection{Convergence of the Coriolis operator}

The convergence of the Coriolis operator was investigated as follows.
A stream function equal to $\cos \varphi \sin \lambda$ was sampled at
dual vertices ($\hat{\Psi}$), enabling exact dual edge normal fluxes
$\widehat{U^\perp} = - \overline{D}_1 \hat{\Psi}$ to be computed.
The same stream function was also sampled at primal vertices ($\Psi$),
enabling exact primal edge normal fluxes $U = - D_1 \Psi$ to be calculated;
approximate dual edge normal fluxes are then given by the
Coriolis operator:
$H \widehat{U^\perp}_\mathrm{approx} = M U^\perp_\mathrm{approx} = -W U = W D_1 \Psi$.
The difference between the two estimates $W D_1 \Psi + H \overline{D}_1 \hat{\Psi}$
gives a measure of the error in the Coriolis operator.

\begin{table}
\centering
\caption{Convergence of the Coriolis operator on hexagonal and cubed sphere grids.}
\label{tab:coriolis}
\begin{tabular}{|r|l|l|r|l|l|}
\hline
Hex    &                 &            & Cube   &                &           \\
Ncells & $L_\infty$ err  &  $L_2$ err & Ncells & $L_\infty$ err & $L_2$ err \\
\hline
    42 & 0.018           & 0.0092     &     54 & 0.0079         & 0.0049    \\
   162 & 0.0049          & 0.0026     &    216 & 0.0055         & 0.0021    \\
   642 & 0.0018          & 0.00066    &    864 & 0.0039         & 0.00092   \\
  2562 & 0.00078         & 0.00017    &   3456 & 0.0022         & 0.00037   \\
 10242 & 0.00036         & 0.000042   &  13824 & 0.0012         & 0.00014   \\
 40962 & 0.00017         & 0.000011   &  55296 & 0.00060        & 0.000050  \\
\hline
\end{tabular}
\end{table}

Values of the error at different resolutions on the two meshes are shown in table~\ref{tab:coriolis}.
On both meshes the $L_\infty$ errors converge at first order. The $L_2$ errors converge
at second order on the hexagonal mesh and between first and second order on the cubed
sphere mesh. This consistency of the Coriolis operator, in contrast to the finite
volume scheme of \citet{thuburn2014}, was one of the primary motivations for investigating
the finite element approach.

\subsection{Solid body rotation}
\label{sbr}

Test case~2 of \citet{williamson1992} tests the ability of models to represent
large-scale steady balanced flow. The exact solution is known, allowing
errors in $\phi$ and $u$ to be computed. The errors on the two meshes
after 5~days are given in table~\ref{tab:sbr}, along with the time steps used
at different resolutions.

\begin{table}
\centering
\caption{Geopotential errors ($\mathrm{m}^2\mathrm{s}^{-2}$) and
velocity errors ($\mathrm{m}\mathrm{s}^{-1}$) for the solid body
rotation test case.}
\label{tab:sbr}
\begin{tabular}{|r|l|l|l|l|l|}
\hline
Ncells & $\Delta t \ \ \mathrm{(s)}$
                     & $L_2(\phi)$ & $L_\infty(\phi)$ & $L_2(u)$   & $L_\infty(u)$ \\
\hline
Hex    &             &             &                  &            &               \\
   642 &   7200      &    19.62    &    ~~43.40       &  0.290     & 0.774         \\
  2562 &   3600      &   ~~8.59    &    ~~14.52       &  0.0940    & 0.217         \\
 10242 &   1800      &   ~~2.27    &   ~~~~4.01       &  0.0244    & 0.0551        \\
 40962 &  ~~900      &   ~~0.584   &   ~~~~1.13       &  0.00609   & 0.0144        \\
\hline
Cube   &             &             &                  &            &               \\
   864 &   7200      &    35.04    &    ~~87.48       &  0.212     & 0.569         \\
  3456 &   3600      &    10.16    &    ~~18.06       &  0.0754    & 0.235         \\
 13824 &   1800      &   ~~2.57    &   ~~~~4.65       &  0.0194    & 0.0692        \\
 55296 &  ~~900      &   ~~0.639   &   ~~~~1.17       &  0.00484   & 0.0257        \\
\hline
\end{tabular}
\end{table}

On the hexagonal mesh the convergence rate is close to second order or better.
On the cubed sphere mesh it is between first and second order for
$L_\infty(u)$ and close to second order for the other error measures.
The errors are considerably smaller than for the finite volume scheme of
\citet{thuburn2014} (their table~6).

\begin{figure}
\begin{centering}
\includegraphics[width=65mm]{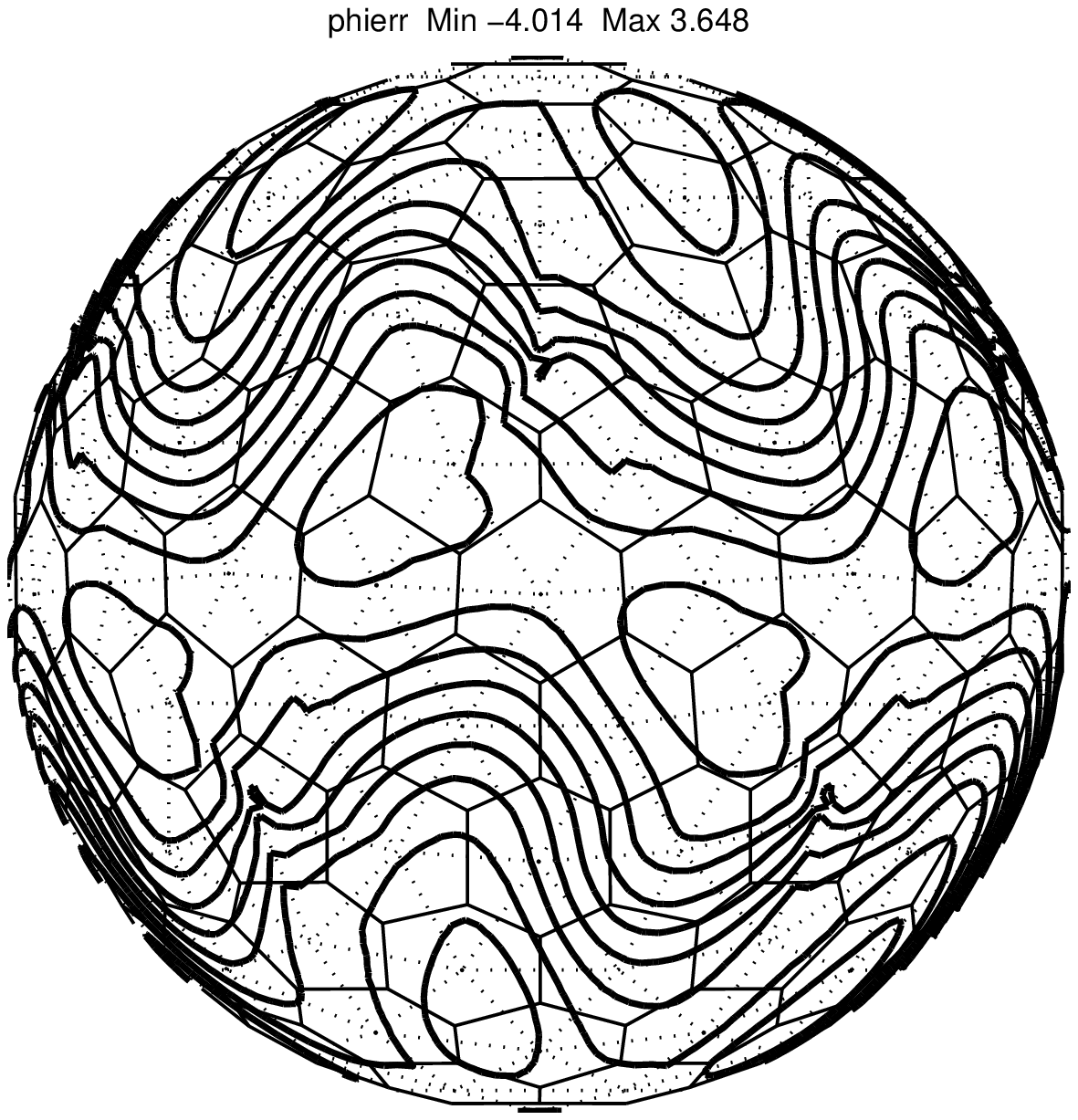}
\includegraphics[width=65mm]{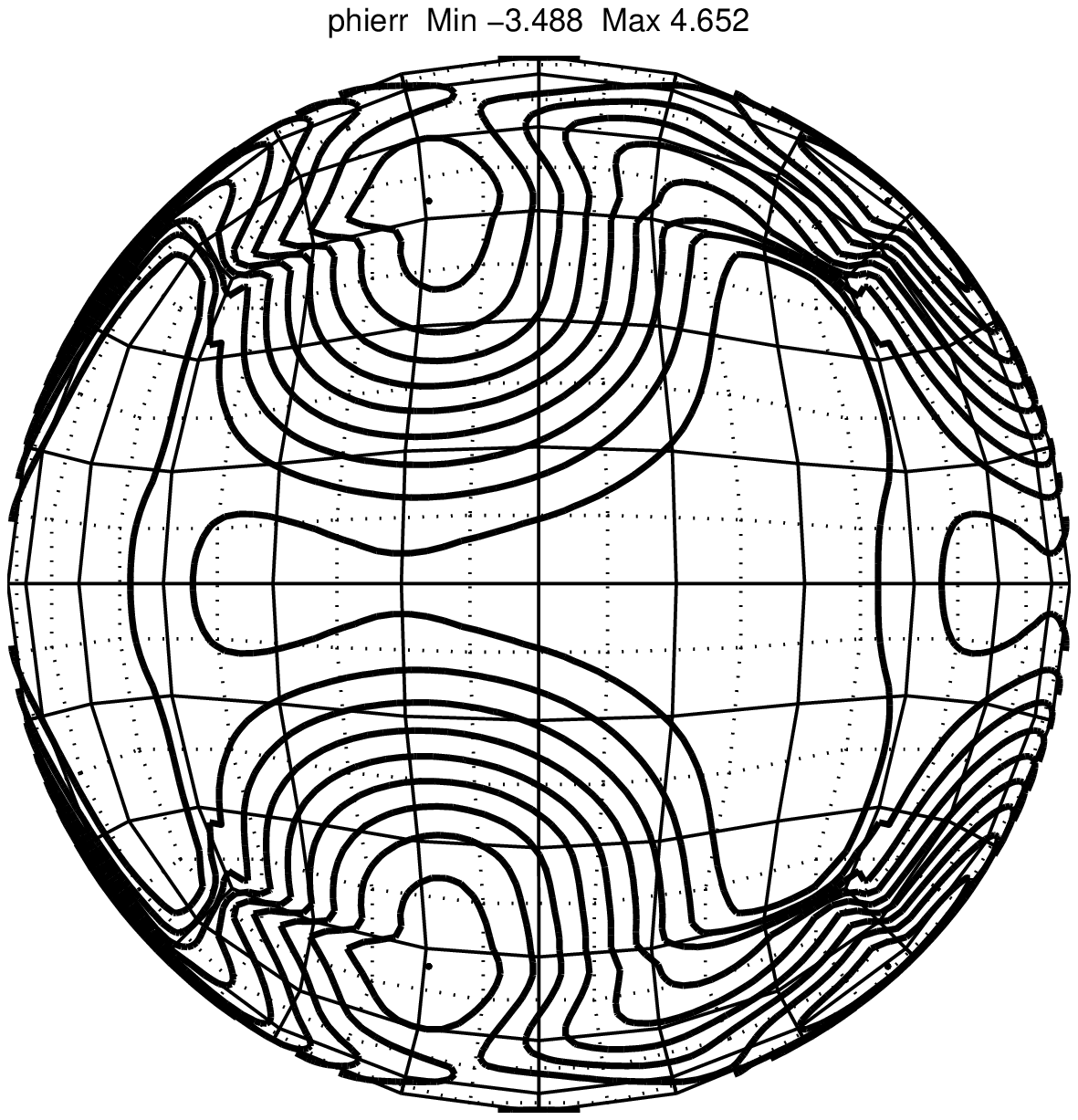}

\caption{Geopotential error ($\mathrm{m}^2\mathrm{s}^{-2}$) after 5~days
for the solid body rotation test case.
Left: hexagonal mesh, 10242~cells. 
Right: cubed sphere mesh, 13824~cells.
In each case 11~evenly space contours (i.e.\ 10 intervals) are used
between the minimum and maximum values. (The coarse resolution meshes
shown as background are for orientation only.)
}
\label{fig:sbr}
\end{centering}
\end{figure}

Figure~\ref{fig:sbr} shows the pattern of geopotential errors after 5~days
at the second highest resolution in the table.
The errors clearly reflect the mesh structure, showing a zonal wavenumber~5
pattern on the hexagonal mesh and a zonal wavenumber~4 pattern on the
cubed sphere mesh. However,
in contrast to the finite volume model, which shows errors concentrated
along certain features of the mesh, the error pattern here is large scale and
almost smooth.

\subsection{Flow over an isolated mountain}
\label{tc5}

Test case~5 of \citet{williamson1992} involves an initial solid body rotation
flow impinging on a conical mid-latitude mountain, leading to the generation of
gravity and Rossby waves and, eventually, a complex nonlinear flow. There is no
analytical solution for this test case, so a high-resolution reference solution
was generated using the semi-implict, semi-Lagrangian ENDGame shallow water
model \citep{zerroukat2009}. The finite element model runs stably with the
time steps given in table~\ref{tab:sbr}, but, as discussed by \citet{thuburn2014}
for the finite volume model and for ENDGame itself, the errors are then dominated
by the semi-implicit treatment of the large amplitude gravity waves present in the
solution. At any given resolution, the errors look almost identical for all combinations
of model and mesh tested. The test was therefore repeated with the time steps
reduced by a factor~4. The resulting height errors at day~15 are shown in table~\ref{tab:tc5}.
The errors on the two meshes are generally very similar, and in most cases are
a little smaller than those produced by the finite volume model
\citep[][table~7]{thuburn2014}.

\begin{table}
\centering
\caption{Height errors ($\mathrm{m}$) for test case~5.}
\label{tab:tc5}
\begin{tabular}{|r|l|l|l|l|l|}
\hline
Ncells & $\Delta t \ \ \mathrm{(s)}$
                     & $L_1(h)$    & $L_2(h)$         & $L_\infty(h)$  \\
\hline
Hex    &             &             &                  &                \\
   642 &   1800      &    36.37    &    ~~50.91       &    191.47      \\
  2562 &  ~~900      &    11.62    &    ~~15.83       &   ~~66.84      \\
 10242 &  ~~450      &   ~~3.12    &   ~~~~4.11       &   ~~15.06      \\
 40962 &  ~~225      &   ~~1.27    &   ~~~~1.82       &  ~~~~9.28      \\
\hline
Cube   &             &             &                  &                \\
   864 &   1800      &    44.11    &    ~~64.93       &    291.35      \\
  3456 &  ~~900      &    17.57    &    ~~25.14       &    100.66      \\
 13824 &  ~~450      &   ~~3.75    &   ~~~~5.25       &   ~~21.42      \\
 55296 &  ~~225      &   ~~1.08    &   ~~~~1.46       &  ~~~~6.47      \\
\hline
\end{tabular}
\end{table}

\begin{figure}
\begin{centering}
\includegraphics[width=100mm]{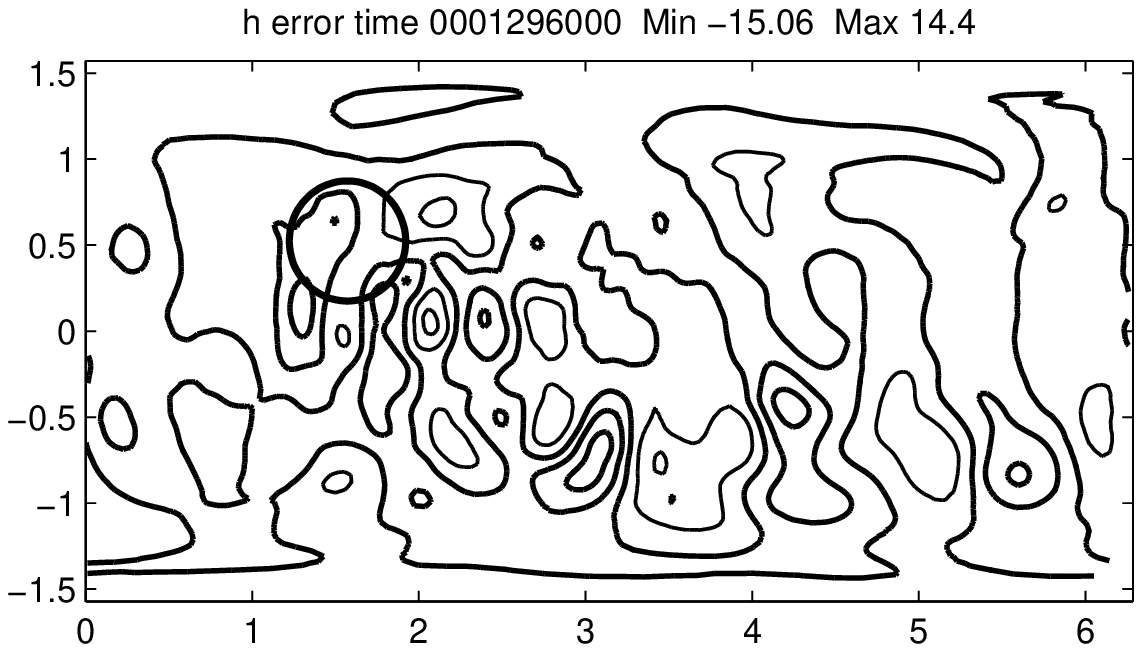}

\includegraphics[width=100mm]{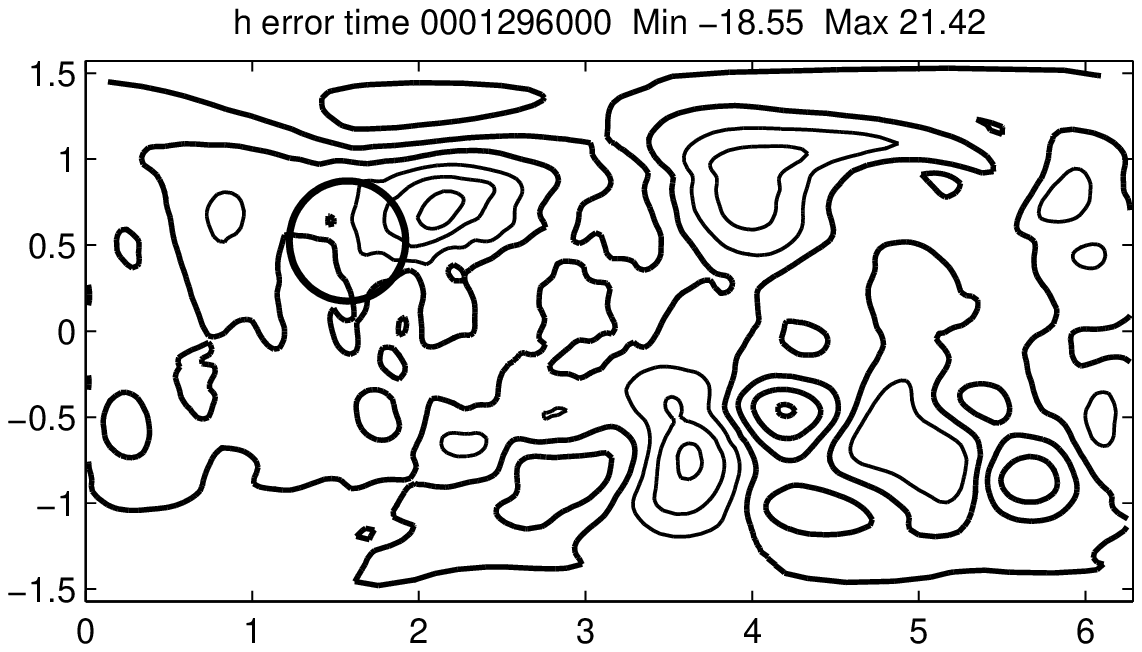}

\includegraphics[width=100mm]{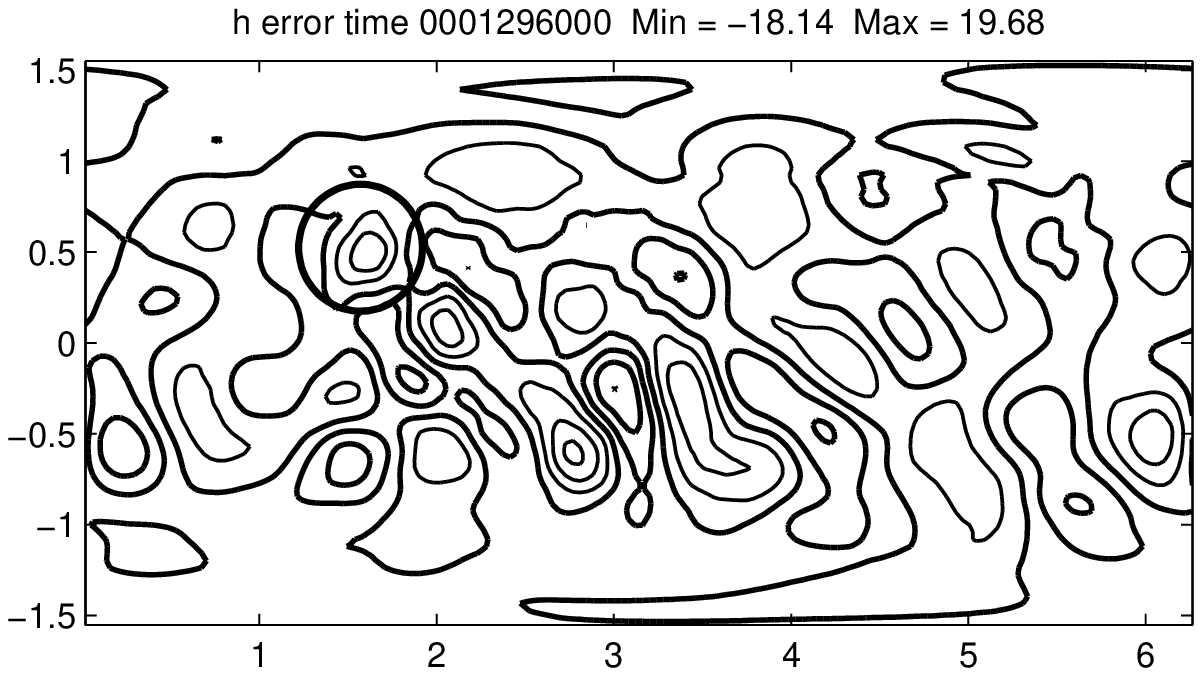}

\caption{Height errors ($\mathrm{m}$) at day~15 for the isolated mountain test case.
Top: hexagonal mesh, 10242~cells. 
Middle: cubed sphere mesh, 13824~cells.
Bottom: ENDGame on a regular longitude-latitude mesh, $160 \times 80$~cells.
The contour interval is $6\,\mathrm{m}$, and zero and negative contours are bold.
The bold circle indicates the position of the mountain.
}
\label{fig:tc5}
\end{centering}
\end{figure}

Maps of height error at day~15 are shown in Fig.~\ref{fig:tc5}.
The errors produced by the finite element model are of comparable
size to those from ENDGame, though the error patterns are different
in the three cases. Comparison with figure~6 of \citet{thuburn2014}
confirms that the errors in the finite element model are somewhat smaller
than those in the finite volume model.

This test case was also run to 50~days at the highest resolutions in table~\ref{tab:tc5}
and several diagnostics relevant to the mimetic properties of the scheme were calculated.
The results are very similar to those shown in figure~8 of \citet{thuburn2014}.
They confirm that mass is conserved to within roundoff error, and that changes in
the total available energy (available potential energy plus kinetic energy) are
much smaller than the conversions between available potential energy and kinetic energy.
The dissipation of available energy and potential enstrophy is associated almost
entirely with the inherent scale-selective dissipation in the advection scheme;
it is very small, of order 1~part per thousand, during the first 20~days, but increases
subsequently as PV contours begin to wrap up and nonlinear cascades become significant,
implying that the inherent dissipation adapts automatically to the flow complexity
in a reasonable way. A dual-mass-like tracer remains consistent with the diagnosed
dual mass field $\bar{\phi}$ to within $2$~parts in $10^4$,
and a PV-like tracer remains consistent with the diagnosed
PV field, to within $3$~parts in $10^3$. The small errors result from
imperfect convergence of various iterative aspects of the solver,
and can be reduced by taking more iterations.

\subsection{Barotropically unstable jet}
\label{galewsky}

The test case proposed by \citet{galewsky2004} follows the evolution of a perturbed
barotropically unstable jet. The case tests the ability of models to handle the
complex small scale vorticity features produced by the rapidly growing instability.
The results are very sensitive to spurious triggering of the instability by
error patterns related to the mesh structure.

\begin{figure}
\begin{centering}
\includegraphics[width=150mm]{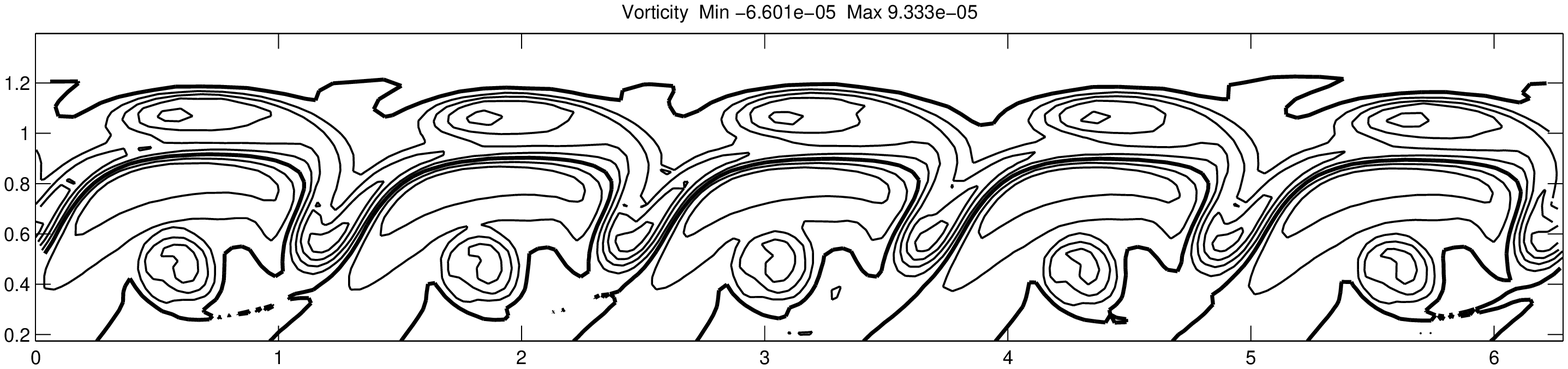}

\includegraphics[width=150mm]{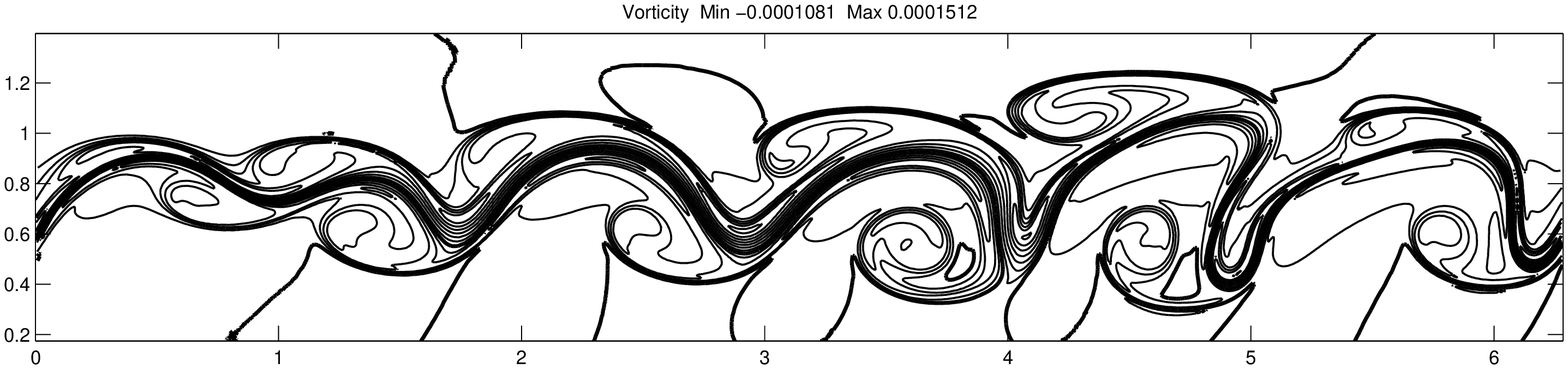}

\includegraphics[width=150mm]{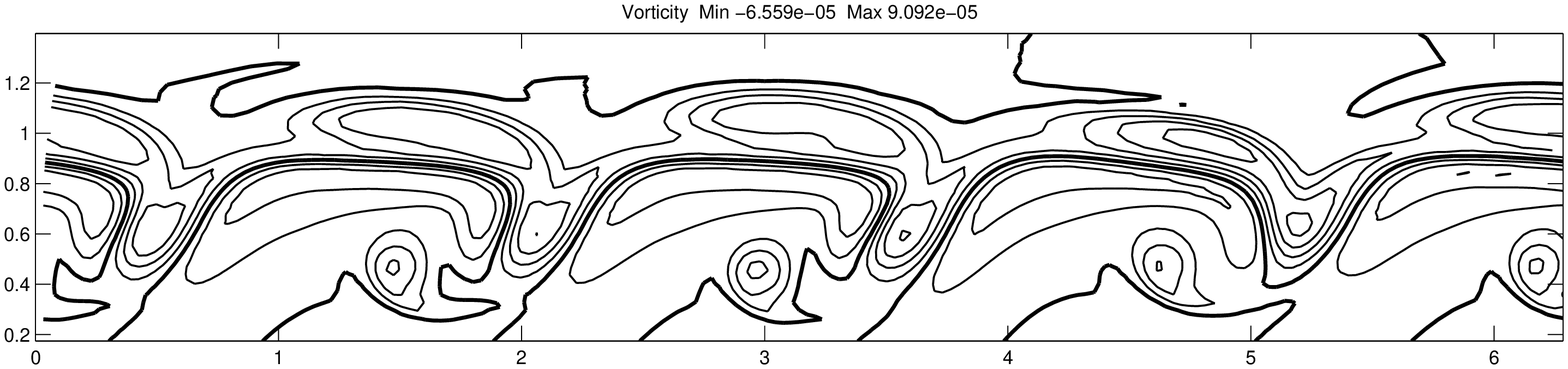}

\includegraphics[width=150mm]{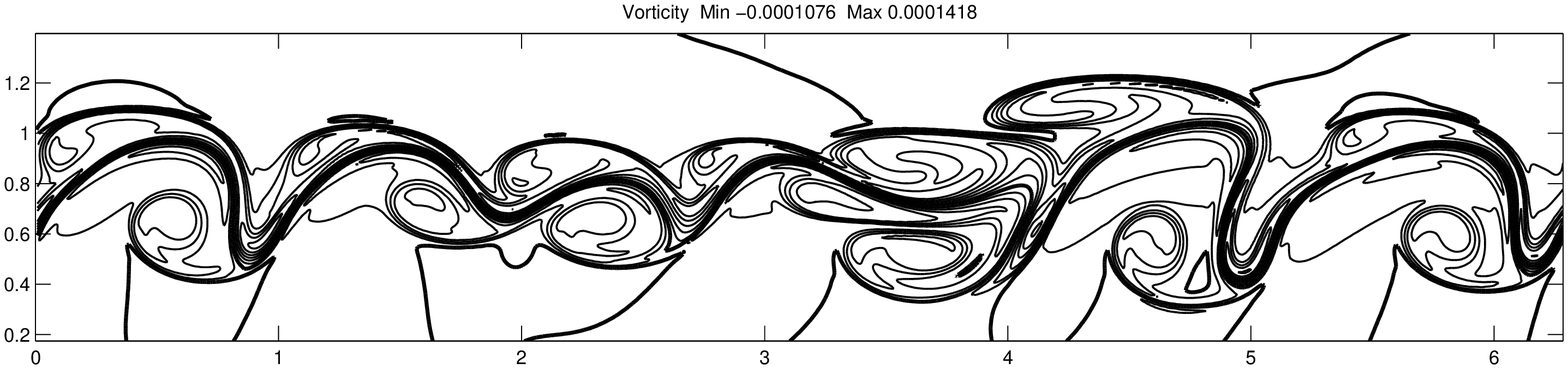}

\includegraphics[width=150mm]{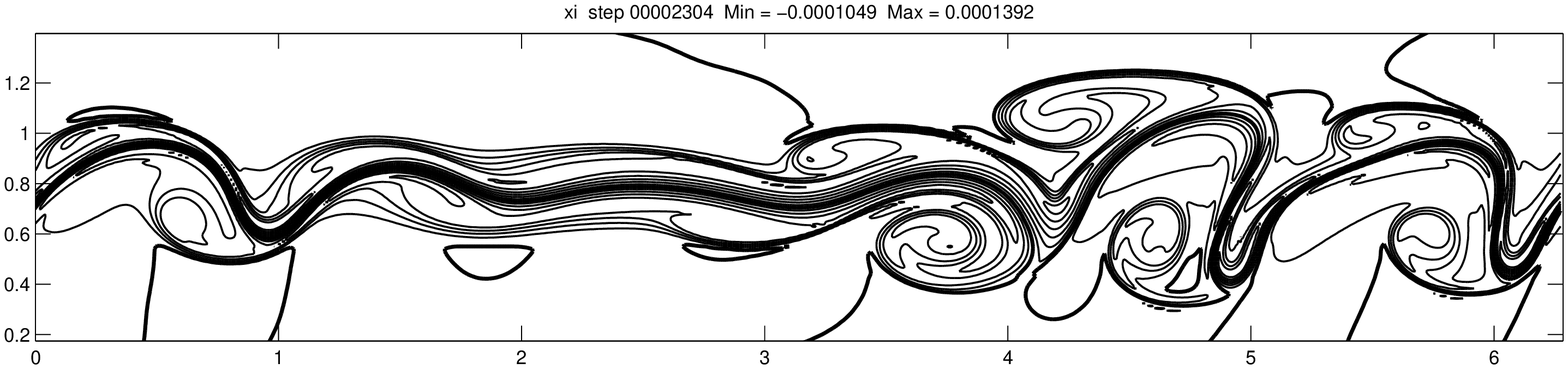}

\caption{Relative vorticity field at day~6 for the barotropic instability test case.
Row~1: hexagonal mesh, 10242~cells, $\Delta t = 900\,\mathrm{s}$.
Row~2: hexagonal mesh, 163842~cells, $\Delta t = 225\,\mathrm{s}$.
Row~3: cubed sphere mesh, 13824~cells, $\Delta t = 900\,\mathrm{s}$.
Row~4: cubed sphere mesh, 221184~cells, $\Delta t = 225\,\mathrm{s}$.
Row~5: ENDGame, $640 \times 320$~cells, $\Delta t = 225\,\mathrm{s}$.
The plotted region is $0^\mathrm{o}$
to $360^\mathrm{o}$ longitude, $10^\mathrm{o}$ to $80^\mathrm{o}$ latitude.
The contour interval is $2 \times 10^{-5}\,\mathrm{s}^{-1}$.
}
\label{fig:galewsky}
\end{centering}
\end{figure}

Figure~\ref{fig:galewsky} shows the relative vorticity field at day~6 for the hexagonal mesh
with 10242~cells and 163842~cells, the cubed sphere mesh with 13824~cells and
221184~cells, and, for comparison, from ENDGame on a $640 \times 320$ longitude-latitude
mesh.
In all cases the vorticity field is free of noise and spurious ripples.
However, at coarse resolution the finite element model solutions show distinct
`grid imprinting', with a
zonal wavenumber~5 pattern on the hexagonal mesh and a zonal wavenumber~4 pattern
on the cubed sphere mesh. At finer resolution the solutions are more similar to the
ENDGame solution, but still show significant development in the longitude range
$\pi/2$~to~$\pi$ where the jet in the ENDGame solution remains quiescent.
The solutions on the hexagonal mesh, especially at the finer resolution, are remarkably
similar to those from the finite volume model \citep[][figure~9]{thuburn2014}.
On the other hand, the solutions on the cubed sphere mesh show some noticable
differences from the finite volume model. At the finer resolution, outside the region
strongly affected by the spurious devleopment, the structure of vorticity features is
slightly more accurate in the finite element model.

\section{Conclusions and discussion}
\label{conclusions}

A method of constructing low-order mimetic finite element spaces
on arbitrary two-dimensional polygonal meshes, using compound elements,
has been presented,
along with corresponding discrete Hodge star operators for mapping
between primal and dual function spaces. The method has been used
as the basis of a numerical model
to solve the shallow water equations on a rotating sphere.
The model has the same mimetic properties, which underpin the
ability to capture important physical properties,
as the finite volume model of \citet{thuburn2014}, but with improved accuracy.

The finite volume model of \citet{thuburn2014} relies on certain properties of the
mesh for accuracy, namely the \citet{heikes1995b} optimization on the hexagonal mesh
and the placement of primal vertices relative to dual vertices on the cubed sphere
mesh. Although identical meshes have been used here to ensure the cleanest possible
comparison, the mimetic finite element scheme does not depend on such mesh properties
for accuracy; thus it provides greater flexibility in the choice of mesh.

An important practical consideration is the computational cost of the method. 
As a rough guide, the cost of the finite element model on a single processor
varied between~3.3 and 4.6~times the cost of ENDGame for the cubed sphere mesh
and between~4.2 and 7.3~times the cost of ENDGame for the hexagonal mesh,
at the resolutions tested\footnote{Martin Schreiber (pers.\ comm.)
reports that the cost of the finite
element model can be significantly reduced, by roughly a factor~2, by reordering the
dimensions of a couple of key arrays to improve cache usage.}.
(For comparison, the cost of the finite volume model
varied between~2.7 and 3.7~times the cost of ENDGame for the cubed sphere mesh
and between~3.3 and 4.9~times the cost of ENDGame for the hexagonal mesh.)
The greater cost on the hexagonal mesh compared to the cubed sphere results from a
combination of a greater stencil size for some operators and, in the
current implementation, a less cache-friendly mesh numbering (the latter could
straightfowardly be optimized).
Given the potential to optimize the implementation and the expected
gains in parallel scalability from the quasi-uniform mesh,
these figures suggest that, despite the need for indirect addressing and the
need to invert several linear operators, the finite element method need
not be prohibitively expensive compared to methods currently used for operational
forecasting, typified by ENDGame.


Computing integrals over compound elements is more complex and costly than
for the usual triangular or quadrilateral elements.
In the current implementation, all the operators
$L$, $M$, $H$, $J$, $R$, $W$, and $T$ are precomputed, thus avoiding the need
for any run-time quadrature in the finite element parts of the
calculations\footnote{Some run-time quadrature is done in the advection scheme
to compute swept area integrals.}.
(Also, once these operators are computed, there is no need to retain the details of
how the compound elements were built from subelements.)
This precomputation is possible because all but one of these operators are linear;
the only nonlinear term (other than advection) is a simple quadratic nonlinearity in
the kinetic energy. In a system with more complex nonlinearities, such as
the pressure gradient term in a compressible three-dimensional fluid, precomputation
might not be possible and some run time quadrature would be unavoidable.
Even so, in a high performance computing environment it is not clear whether
precomputation or run-time quadrature would be most efficient, given the
relative cost of memory access and computation (David Ham, pers.\ comm.).

The mathematical similarity of the finite element and finite volume
formulations has been emphasized, the principal difference being the
appearance of mass matrices in the finite element formulation.
The similarity is made even clearer if we use \refeqn{eqn_hodge02},
\refeqn{eqn_hodge11} and \refeqn{DIeqHD} to rewrite \refeqn{mudot}
in the equivalent dual space form
\begin{equation}
\label{mudotdual}
\dot{\hat{U}} + \widehat{Q^\perp} + \overline{D}_1 ( \hat{\Phi} + \hat{K} )  =  0 .
\end{equation}
The velocity degrees of freedom then correspond to dual edge circulations,
and $\hat{U}$ can be identified with the $V$ of \citet{thuburn2014}.
Equations \refeqn{phidot} and \refeqn{mudotdual} explicitly involve only
the topological derivative operators $D_2$ and $\overline{D}_1$;
the metric enters through the Hodge star operators needed to map between
primal and dual function spaces. This approach of isolating the metric
from the purely topological operators in order to construct numerical
methods with mimetic properties on complex geometries or meshes
has been advocated by several authors
\citep[e.g.][and references therein]{bossavit1998,hiptmair2001,palha2014}.

It is also worth emphasizing that the roles of primal and dual function
spaces are not symmetrical here. Although any given field may be represented
in both the primal and dual spaces, with a reversible Hodge star
map between them, only primal space test functions are ever used, and
so only primal space mass matrices appear, and dual space weak derivatives
are never needed. (An interesting alternative would be to use only dual
space test functions; then the prognostic equations remain
\refeqn{phidot} and \refeqn{mudotdual}, but \refeqn{eqn_hodge02}-\refeqn{eqn_hodge20}
are replaced by
\begin{equation}
I^T \Phi = \hat{L} \hat{\Phi} ,
\label{eqn_alt_hodge02}
\end{equation}
\begin{equation}
H^T U = \hat{M} \hat{U} ,
\label{eqn_alt_hodge11}
\end{equation}
\begin{equation}
J^T \Xi = \hat{N} \hat{\Xi} ,
\label{eqn_alt_hodge20}
\end{equation}
where $\hat{L}$, $\hat{M}$ and $\hat{N}$ are the mass matrices for the spaces
$\mathbf{V}^0$, $\mathbf{V}^1$ and $\mathbf{V}^2$, respectively.)

Only the lowest order polygonal finite element spaces are used here:
compound $\text{P1-RT0-P0}^\text{DG}$. An interesting question is whether
the approach can be extended to higher order. 
The harmonic extension idea of \citet{christiansen2008} has been extended
to higher order by \citet{christiansen2010}.
It appears plausible that
higher order compound elements could be built from constrained linear
combinations of, for example, the $\text{P2}^+\text{-BDFM1-P1}^\text{DG}$
elements recommended by \citet{cotter2012}, but the details have yet
to be worked out.
A more subtle question is whether suitable higher order dual spaces
can be constructed.

Another, more straightforward, extension of the compound element approach
is to three dimensions. The compound elements used here can be extruded
into polygonal prisms; we have made some initial progress in working out
the details of using such a scheme for the compressible Euler equations.
(In atmosphere and ocean models it is desirable,
for several reasons, to use a columnar mesh.)
Fully three-dimensional compound elements can also be constructed using
the discrete harmonic extension approach. These might be useful,
for example, to implement a finite element version of the cut cell
method for handling bottom topography
\citep[e.g.][and references therein]{lock2012}
while retaining a columnar mesh.

Besides their ability to use arbitrary polygonal meshes, another potentially
useful property of the compound elements used here is that the function
spaces are built directly in physical space, without the need for
Piola transforms. Thus, for example, globally constant functions
are always contained in $\mathbf{V}_2$.
In this way, the compound elements avoid the reduced convergence rate, and
even loss of consistency, discussed by \citet{arnold2014}, and so provide an
alternative to the {\it rehabilitation} technique of \citet{bochev2008}.

\vspace{3mm}

{\bf Acknowledgements}

We are grateful to Thomas Dubos for drawing our attention to the work of
\citet{buffa2007} and \citet{christiansen2008}, to Martin Schreiber
and David Ham for valuable discussions on code optimization,
and to Nigel Wood for helpful comments on a draft of this paper.
This work was funded by the
Natural Environment Research Council under the ``Gung Ho'' project
(grants NE/I021136/1, NE/I02013X/1, NE/K006762/1 and NE/K006789/1).

\appendix

\section{Operator inverses and sparse approximate $M$ inverse}
\label{lumped}

Inverses of the $H$ and $J$ operators are needed at the beginning of
every time step, and inverses of $H$ and $M$ are needed at every Newton
iteration. These are computed by (under- or over-relaxed) Jacobi
iteration based on a diagonal approximation of the relevant operator.
E.g., to solve $A x = R$, define
\begin{equation}
\label{jacobi_fg}
x^{(1)} = (A^{*})^{-1} R
\end{equation}
where $A^{*}$ is a diagonal approximation to $A$, then iterate:
\begin{equation}
x^{(l+1)} = x^{(l)} + \mu (A^{*})^{-1} (R - A x^{(l)}).
\end{equation}

A diagonal approximation $J^{*}$ to the operator $J$ is defined by
demanding that, for every dual cell $j$, $J^{*}$ and $J$ should give the
same result in dual cell $j$ when acting on the $\mathbf{V}^2$
representation of a constant scalar field.
A diagonal approximation $M^{*}$ to the velocity mass matrix $M$
is defined by demanding that, for every edge $e$, $M^{*}$ and $M$ should give
the same result at edge $e$ when acting on the $\mathbf{V}_1$ representation
of a solid body rotation velocity field whose maximum velocity
is normal to primal edge $e$.
A diagonal approximation $H^{*}$ to the operator $H$
is defined by demanding that, for every edge $e$, $H^{*}$ and $H$ should give
the same result at edge $e$ when acting on the $\mathbf{V}^1$ representation of
a solid body rotation velocity field whose maximum velocity is
tangential to dual edge $e$.

Optimal values of the relaxation parameter $\mu$ were found to depend
on the operator and mesh structure. The values used are given in table~\ref{tab:relax}.
\begin{table}
\centering
\caption{Relaxation parameters used for Jacobi iteration for
operator inverses.}
\label{tab:relax}
\begin{tabular}{|l|c|c|c|}
\hline
Grid  & $J^{-1}$           & $M^{-1}$    & $H^{-1}$    \\
\hline
Hex   &  1.4               &  1.4        &  1.4       \\
Cube  &  1.4               &  0.9        &  1.4       \\
\hline
\end{tabular}
\end{table}
For the inverses that occur once per time step, 10 Jacobi iterations are used.
For those that occur at every Newton iteration, 2 Jacobi iterations are used
taking the solution obtained at the previous Newton iteration as the first guess
(or \refeqn{jacobi_fg} on the first Newton iteration).

A sparse approximate inverse $\mathcal{M}^{-1}$ of the $\mathbf{V}_1$ mass matrix
is needed for the Helmholtz problem. On the hexagonal mesh
it is sufficient to use a diagonal approximation
\begin{equation}
\mathcal{M}^{-1} = (M^{*})^{-1}.
\end{equation}
However, on the cubed sphere mesh, whose dual and primal edges are not
mutually orthogonal, such a diagonal approximation is less accurate
and limits the convergence of the Newton iterations.
Therefore we use instead an approximation based on a single Jacobi iteration
towards the inverse of $M$:
\begin{equation}
\mathcal{M}^{-1} = (M^{*})^{-1} \left\{ (1 + \mu)\mathrm{Id} - \mu M (M^{*})^{-1} \right\},
\end{equation}
where $\mathrm{Id}$ is the identity matrix.
This approximate inverse is not diagonal but has the same stencil as
$M$ itself.

\vspace{10mm}
\noindent
{\bf References}
\vspace{4mm}



\end{document}